\documentclass[reqno, 11pt]{amsart}

\usepackage{amssymb,amscd,bm}
\usepackage[mathscr]{eucal}
\usepackage[dvips]{graphicx,color}
\input psfig.sty

\numberwithin{equation}{section}

\newtheorem{theorem}{Theorem}[section] \newtheorem{lemma}[theorem]{Lemma}
\newtheorem{definition}[theorem]{Definition}

\newtheorem{proposition}[theorem]{Proposition}

\newtheorem{corollary}[theorem]{Corollary}

\newcommand{\ZZ}{{\mathbb Z}}
\newcommand{\RR}{{\mathbb R}}
\newcommand{\TT}{{\mathbb T}}

\newcommand{\NN}{{\mathbb N}}

\newcommand{\pr}{{\mathsf P}}
\newcommand{\dif}{\mathrm{d}}

\newcommand{\UU}{{\mathcal U}}
\newcommand{\Sc}{{\mathcal S}}
\newcommand{\HH}{{\mathscr H}}
\newcommand{\MM}{{\mathscr M}}
\newcommand{\LL}{{\mathscr L}}
\newcommand{\EE}{{\mathscr E}}
\newcommand{\Ts}{{\mathscr T}}
\newcommand{\As}{{\mathscr A}}
\newcommand{\Leb}{\mathfrak{m}}
\newcommand{\vu}{\nu}

\newcommand{\sd}{$\Sigma\Delta$ }

\begin{document}

\title[Ergodic Dynamics in \sd Quantization]
{Ergodic Dynamics in \sd Quantization:\\
Tiling Invariant Sets and Spectral Analysis of Error}

\date{August 20, 2003}

\author{Nguyen T.~Thao and C.~S{\.\i}nan G\"unt\"urk
}

\address{Department of Electrical Engineering, City
College and Graduate School, City University of New York,
Convent Avenue at 138th Street, New York, NY 10031}
\email{thao@ee-mail.engr.ccny.cuny.edu}

\address{Courant Institute of Mathematical Sciences,
New York University, 251 Mercer Street, New York, NY 10012.}
\email{gunturk@cims.nyu.edu}

\thanks{This work has been supported in part by the
National Science Foundation Grants DMS-0219072, DMS-0219053
and CCR-0209431.}

\begin{abstract}
This paper has two themes that are intertwined: The first is
the dynamics of certain piecewise affine maps on $\RR^m$
that arise from a class of analog-to-digital conversion methods called
\sd quantization. The second is the analysis of reconstruction error
associated to each such method.
\par
\sd quantization generates approximate representations of
functions by sequences that lie in a restricted set of discrete
values. These are special sequences in that their local averages
track the function values closely, thus enabling simple
convolutional reconstruction. In this paper, we are concerned with
the approximation of constant functions only, a basic case that
presents surprisingly complex behavior. An $m$th order \sd scheme
with input $x$ can be translated into a dynamical system that
produces a discrete-valued sequence (in particular, a $0$--$1$
sequence) $q$ as its output. When the schemes are stable, we
show that the underlying piecewise affine maps possess
invariant sets that tile $\RR^m$ up to a finite multiplicity.
When this multiplicity is one (the single-tile case), the dynamics
within the tile is isomorphic to that of a generalized skew
translation on $\TT^m$.
\par
The value of $x$ can be approximated using any consecutive $M$
elements in $q$ with increasing accuracy in $M$. We show that the
asymptotical behavior of reconstruction error depends on the
regularity of the invariant sets, the order $m$, and some
arithmetic properties of $x$. We determine the behavior in a
number of cases of practical interest and provide good upper
bounds in some other cases when exact analysis is not yet
available.
\end{abstract}

\maketitle

%\newpage
\setcounter{tocdepth}{1}
%\tableofcontents

%\pagenumbering{roman}
%\tableofcontents \newpage
%\pagenumbering{arabic}

\section{Introduction}
\label{Intro}

\par
This paper is motivated by the mathematical problems exhibited in and
suggested by a class of real-world practical algorithms that are used
to perform analog-to-digital conversion of signals. There will be
two themes in our study of these mathematical problems. The first
theme is the dynamics of certain piecewise affine maps on $\RR^m$
that are associated with these algorithms.
The second theme is the analysis of the reconstruction error.
While the first theme is somewhat independent of the second and is
of great interest on its own, the second
theme turns out to be crucially dependent on the first and is
of interest for theoretical as well as practical reasons.

\par
Let us start with the following abstract
algorithm for analog-to-digital encoding:
For each input real number $x$ in some interval $I$,
there is a map $\Ts_x$ on a space $\Sc$,
and a finite partition $\Pi_x = \{\Omega_{x,1},\dots,\Omega_{x,K}\}$ of $\Sc$.
For a fixed set of real numbers
$d_1 < \cdots < d_K$, and a typically fixed (but arbitrary)
initial point $u_0 \in \Sc$,
we define a discrete-valued output sequence $q := q_x$ via
\begin{equation}
q[n] = d_i \;\;\;\; \mbox{ if } \;\;\;\; u[n{-}1]:=\Ts_x^{n-1}(u_0)
\in \Omega_{x,i}.
\end{equation}
We would like the mapping $x \mapsto q$ to be invertible in
a very special way: For an {\sl input-independent} family of averaging
convolutional
kernels $\phi_M \in \ell^1(\ZZ)$, $M=1,2,\dots$, we require that
for all $x \in I$, as $M \to \infty$
\begin{equation}\label{abs-approx}
(q * \phi_M)[n] \longrightarrow x, \;\;\;\;\;\mbox{ uniformly in } n.
\end{equation}
For normalization, we ask the size of the averaging window (the
support of $\phi_M$) to grow linearly in $M$,\footnote{ It will be
of interest to use infinitely supported kernels as well. We will
define the necessary modifications to handle this situation
later.} and the weights to satisfy $\sum_n \phi_M[n] = 1$.
\par
Note that such an encoding of real numbers is inherently different
from binary-expansion (or any other expansion in a number system)
in that, due to (\ref{abs-approx}), equal length
portions of the sequence $q$ are required to be
equally good in approximating the value
of $x$. Hence, there is a ``translation-invariance'' property
in the representation.
\par
This setting is a special case of a more general one in which
$x = (x[n])_{n\in\ZZ}$ is a bounded sequence
taking values in $I$ and
\begin{equation}
q[n] = d_i \;\;\;\; \mbox{ if } \;\;\;\;
u[n{-}1] \in \Omega_{x[n],i},
\end{equation}
where we now define $u[n] := \Ts_{x[n]}(u[n{-}1])$, and require that
\begin{equation}\label{general-abs-approx}
(q-x)* \phi_M \longrightarrow 0 \;\;\;\;\;\; \mbox{ uniformly.}
\end{equation}
The basic motivation behind this type of encoding is the following
intuitive idea: Let the elements $x[n]$ be closely and regularly spaced
samples of a smooth function $X:\RR \to I$. Since
local averages of these samples around any point $k$ would approximate $x[k]$,
i.e., $x*\phi_M \approx x$ for suitable $\phi_M$,
(\ref{general-abs-approx}) would then
imply that the sequence $x$ (and therefore the function $X$) can
be approximated by the convolution $q*\phi_M$.

\par
Such analog-to-digital encoding algorithms have been developed and
used in electrical engineering for a few decades now. Most notable
examples are the {\em \sd quantization} (also called \sd {\em modulation})
of audio signals and
the closely related {\em error-diffusion} in digital halftoning of images.
There are several sources in the electrical engineering literature
on the theoretical and practical aspects of \sd quantization
\cite{sigmadelta1, gray1, sigmadelta2}. Digital
halftoning and its connections to \sd quantization can be found in
\cite{Adler, Anastas, Bernard, Kite1, Ulichney}. Recently, \sd quantization
has also received interest in the mathematical community,
especially in approximation theory and information theory, since
a very important question is the rate of convergence in
(\ref{general-abs-approx}) \cite{ingrid_democ, ingrid_devore, Sinan1, 
exp_decay}.

\par
We give in Section \ref{describe-sd}
the original description of an $m$th order \sd modulation scheme 
in terms of difference equations. The underlying specific map $\Ts_x$,
which we then refer to as $\MM_x$ (the ``modulator map'')
is described in Section \ref{piecewise-affine}; $\MM_x$ is the
piecewise affine transformation on $\Sc = \RR^m$ defined by
\begin{equation}\label{Mx}
\MM_x({\bf v}) = 
%\As_{x,i}({\bf v}) := 
{\bf L}{\bf v} + (x-d_i){\bf 1} \;\;\;\;
\mbox{ if } \;\; {\bf v} \in \Omega_{x,i},
\end{equation}
where ${\bf L}:={\bf L}_m$ is the $m{\times}m$ lower triangular
matrix of $1$'s and ${\bf 1} := {\bf 1}_m := (1,\dots,1) \in
\ZZ^m$. Each \sd scheme is therefore 
characterized by its order $m$, the partition $\Pi_x$, and the 
numbers $\{d_i\}$.
A scheme is called $k$-bit if the size $K$ of the partition $\Pi_x$
satisfies $2^{k-1} < K \leq 2^k$. If the numbers $\{d_i\}$ are in an
arithmetic progression, this is referred to as {\em uniform
quantization}. As a consequence of the normalization $\sum_n
\phi_M[n] = 1$, the input numbers $x$ are chosen in $I \subset
[d_1,d_K]$. A scheme is said to be {\em stable} if for each $x$,
forward trajectories under the action of $\MM_x$ are bounded in
$\RR^m$. (More refined definitions of stability will be given in Section
\ref{piecewise-affine}.) The partition $\Pi_x$ is an essential part of the
algorithm for its central role in stability.

\par
It is natural to measure the accuracy of a scheme by how fast the
worst case error $\| (x - q) * \phi_M \|_\infty$
converges to zero. It is known
that for an $m$th order stable scheme, and an appropriate choice for the
family $\mathscr{F}= \{\phi_M\}$ of filters,\footnote{We shall adopt the
electrical engineering terminology ``filter''
to refer to a sequence (or function) that acts convolutionally.}\,
this quantity is $O(M^{-m})$ \cite{ingrid_devore}.
The hidden constant depends on the scheme
as well as the input sequence $x$. Here,
the exponent $m$ is not sharp; in fact,
for $m=1$ and $m=2$, improvements have been given for various
schemes \cite{Sinan1, SinanThao1}. We will review the basic approximation
properties of \sd quantization in Section \ref{describe-sd}.

\par
In applications, it is also common to measure the error in the root mean
square norm due its more robust nature 
(this norm is defined in Section \ref{sec-error}). 
It is known for a small class of schemes we call {\em
ideal}, and a small class of sequences (basically, constants and
pure sinusoids) that this norm, when averaged over a smooth distribution
of values of $x$, has the asymptotic behavior $O(M^{-m-1/2})$
\cite{GrayMulti,gray2,He92}. The analyses employed in obtaining
these results rely on very special properties of the ideal
schemes, such as employing an (effective) $m$-bit uniform quantizer
for the $m$th order scheme.
It was not known how to extend these results to
low-bit schemes (in particular, $1$-bit schemes) of high order for which
experimental results and simulation suggested similar asymptotical
behavior for the root mean square error.

\par
It is the topic of this paper to provide a general framework and
methodology to analyze \sd quantization in an arbitrary setup
(in terms of partition and number of bits)
when inputs are constant sequences.
With regards to the first theme of this paper, we prove in
Section \ref{sec-tiling} that the maps $\MM_x$ have an
outstanding property of yielding {\em tiling invariant sets},
up to a multiplicity that is determined by the map.
In the particular case of single tiles being invariant under $\MM_x$
(which also appears to be systematically satisfied by all practical \sd 
quantization schemes), we
develop a spectral theory of \sd quantization. This constitutes 
the second theme of the paper. The
particular consequence of tiling that enables our spectral
analysis is presented in Section \ref{conseq-tiling}. The
resulting new error analysis for general and particular cases is
presented in the remainder of the paper.

\subsection*{Some notation}
\par The symbols $\RR$, $\ZZ$, and $\NN$
denote the set of real numbers, the set of integers and the set
of natural numbers, respectively.
$\TT$ denotes the set of real numbers modulo $1$, i.e., $\TT = \RR/\ZZ$.
Functions on $\RR^m$ that are $1$-periodic in each dimension
are  assumed to be
defined on $\TT^m$ via the identification $\TT=[0,1)$, and
functions defined on $\TT^m$ are extended to $\RR^m$ by periodization.
\par Vectors and matrices are denoted in boldface letters.
Transpose is denoted by an upperscript $\top$.
The $j$'th coordinate of a vector ${\bf v}$
is denoted by $v_j$, unless otherwise specified.
Sequence elements are denoted using brackets, such as
in $\omega = \big ( \omega[n] \big)_{n\in\ZZ}$.
The sequence $\tilde{\omega}$ denotes time reversal of $\omega$
defined by $\tilde \omega[n]:=\omega[-n]$, and
the symbol $*$ is used to denote the convolution operation.
\par We define two types of autocorrelation. The
autocorrelation $A_f$ is defined
for square integrable functions (or sequences) $f$, by the formula
$$A_f(t) = (f*\tilde f)(t) = \int f(\xi)\overline{f(\xi+t)}\,\dif \xi.$$
The autocorrelation $\rho_{\omega}$ is defined for
bounded sequences (or functions) $\omega$, by the formula
$$ \rho_{\omega}[k] = \lim_{N\rightarrow\infty}~
\frac{1}{N}\sum_{n=1}^{N} \omega[n]\overline{\omega[n+k]},$$
provided the limit exists.
\par The Fourier series coefficients of a measure $\mu$ on $\TT$
are given by
$$ \hat \mu[n] := \int_\TT e^{-2\pi i n\xi}\, \dif \mu(\xi).$$
Whenever convenient,
the Fourier transform of a sequence $s = (s[n])_{n \in \ZZ}$ will be denoted
by the capital letter $S$, i.e., $\hat S = s$.
\par The ``big oh'' $f=O(g)$ and the ``small oh'' $f=o(g)$ notations will
have their usual meanings. When constants matter, we also use
the notation $f \lesssim_\alpha g$ to denote that there
exists a constant $C$ that may depend on the parameter (or set of
parameters) $\alpha$ such that $f \leq Cg$. We write $f \asymp g$ if
$f \lesssim g$ and $g \lesssim f$, which is the same as $f = \Theta(g)$.

\section{Basic theory of \sd quantization}
\label{describe-sd}

In this section, we describe the principles of \sd quantization
(modulation)
via a set of defining difference equations.
The description in terms of piecewise affine maps on $\RR^m$
will be given in Section \ref{piecewise-affine}.
Although the schemes representable by these difference
equations do not constitute the whole collection of algorithms
called by the name \sd modulation, they are sufficiently general
to cover a large class of algorithms that are used in practice and many
more to be investigated.
\par
\begin{figure}[t]
\begin{center}
\includegraphics[height = 1.8in]{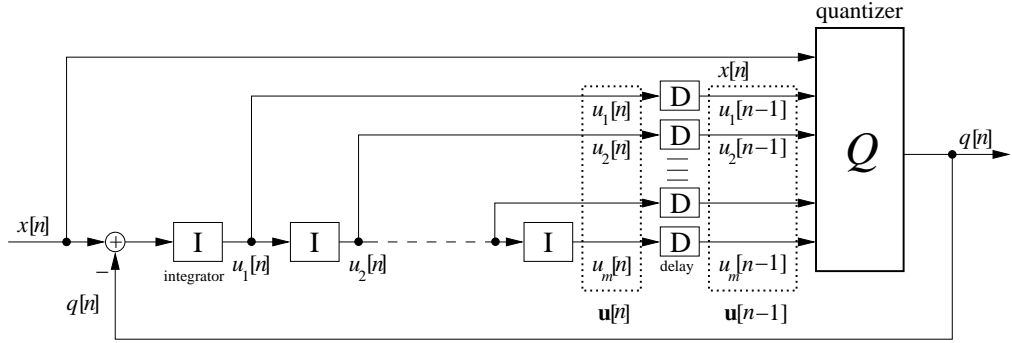}
\end{center}
\caption{\label{sdequiv} Block diagram of an $m$th order \sd modulator.}
\end{figure}

Let $m$ be the order of the scheme, and
$x = (x[n])_{n\in\ZZ}$ be the input sequence. Then a sequence
of state-vectors, denoted
$$ {\bf u}[n] = \Big ( u_1[n], \,\dots,\, u_m[n] \Big )^\top, \;\;\;\;
n=0,1,\dots$$
and a sequence of output quantized values (or symbols), denoted $q[n]$,
$n=1,2,\dots$,
are defined recursively via the set of equations
\begin{equation}\label{state-eq-1}
\left\{
\begin{array}{rcll}
q[n] & = & Q(x[n],{\bf u}[n{-}1]), & \\
u_1[n] & = & u_1[n{-}1] + x[n]-q[n] & \\
u_j[n] & = & u_j[n{-}1] + u_{j{-}1}[n], & 2 \leq j \leq m,
\end{array}\right.
\end{equation}
where the mapping $Q: \RR^{m+1} \to \{d_1,\dots,d_K\}$, called
the {\em quantization rule}, or simply the {\em quantizer} of the
\sd modulator, is specific to the scheme. In circuit theory, these equations
are represented as a feedback-loop system via the block diagram
given in Figure \ref{sdequiv}.

\par
In addition to producing the output sequence $q$, 
the role of the quantizer $Q$ of a \sd modulator is to keep the
variables $u_j$ bounded. A more precise definition of this notion of
stability will be given later.
Let us see how boundedness of $u_j$ results in a simple
reconstruction algorithm.
It can be seen directly from (\ref{state-eq-1}) that for each $j=1,\dots,m$,
the state variable $u_j$ satisfies
\begin{equation}\label{x-q=du}
x - q = \Delta^j u_j,
\end{equation}
where $\Delta$ is the difference operator
defined by $(\Delta v)[n] = v[n]-v[n{-}1]$. Consider $j=1$, and assume
that $x$ is constant. From this, it follows that
\begin{eqnarray}
\left |x - \frac{1}{M}\sum_{k=n+1}^{n+M} q[k] \right | & = &
\frac{1}{M} \left |\sum_{k=n+1}^{n+M} (x - q[k]) \right | \nonumber \\
& = &
\frac{1}{M} \left |\sum_{k=n+1}^{n+M} (u_1[k]-u_1[k{-}1]) \right | \nonumber \\
& = &
\frac{1}{M} \big |u_1[n{+}M] - u_1[n] \big | \nonumber \\
& \leq &
\frac{2}{M} \big \|u_1 \big \|_\infty.
\end{eqnarray}
This means that simple averaging of {\sl any} $M$ consecutive output
values $q[k]$ yields a reconstruction within $O(M^{-1})$.
\par
This approximation result can
be generalized easily. For simplicity of the discussion, let us assume
that the difference equation
(\ref{x-q=du}) is satisfied on the whole of $\ZZ$ (with some care,
this can be achieved via backwards iteration of  (\ref{state-eq-1})).
For a given averaging filter $\phi \in \ell^1(\ZZ)$ with $\sum\limits_n
\phi[n] = 1$, let
\begin{equation}
e_{x,\phi} := x - q*\phi
\end{equation}
be the error sequence. Since $x$ is a constant sequence, we have
$x = x*\phi$. Therefore
\begin{equation}
\label{error-eq}
e_{x,\phi} = (x - q)*\phi = (\Delta^m u_m)*\phi = u_m * (\Delta^m \phi),
\end{equation}
where at the last step we have used commutativity of convolutional operators.
From this, we obtain
\begin{equation}\label{err-eq-bnd}
\big \| e_{x,\phi} \big \|_\infty \leq \big \| u_m \big \|_\infty
\big \| \Delta^m \phi \big \|_1.
\end{equation}
It is not hard to show that
there is a family of averaging kernels $\phi^{(m)}_M$
(which can be, for instance, discrete B-splines of degree $m$)
with support size growing linearly in $M$
such that $ \| \Delta^m \phi^{(m)}_M \|_1 \leq C_m  M^{-m}$.
Combined with (\ref{err-eq-bnd}),
this yields the bound $O(M^{-m})$ on the uniform approximation error.
A proof of this result
in the more general setting of oversampling of bandlimited
functions can be found in \cite{ingrid_devore,CSG_thesis}.

\section{Mean square error and its spectral representation}
\label{sec-error}

For the rest of this paper, we shall be interested in the mean
square error (also called, the {\em time-averaged square error})
of approximation defined by
\begin{equation}
\EE(x,\phi) :=
\lim_{N\to\infty} \; \frac{1}{N} \sum_{n=1}^N \big | e_{x,\phi}[n] \big |^2,
\end{equation}
provided the limit exists (otherwise the $\lim$ is replaced by a
$\limsup$). The {\em root} mean square error is defined
to be $\sqrt{\EE(x,\phi)}$. For convenience in the notation,
we shall work with $\EE(x,\phi)$.

\par
The mean square error enjoys properties that are desirable from
an analytic point of view. The definition of autocorrelation sequence
yields an alternative description given by
\begin{equation}\label{MSE-exp-0}
\EE(x,\phi) = \rho_{e_{x,\phi}}[0].
\end{equation}
Using the formula (\ref{error-eq}) and the standard relation
$\rho_{\omega*g} = \rho_\omega * g * \tilde g\;$
whenever $\rho_\omega$ exists and $g \in l^1$, we find that
\begin{equation}\label{MSE-exp}
\EE(x,\phi) = (\rho_{u_j}*(\Delta^j \phi)*
\widetilde {(\Delta^j \phi)})[0].
\end{equation}
This formula is valid for any $j=1,\dots,m$, provided
$\rho_{u_j}$ exists. In fact, it suffices to compute $\rho_{u_m}$
only, since $u_j = \Delta^{m-j} u_m$ yields
$$ \rho_{u_j} = \Delta^{m-j} * \widetilde{\Delta^{m-j}}* \rho_{u_m}.$$
We shall abbreviate $\rho_{u_m}$ by $\rho_u$.

\par
The computation of $\EE(x,\phi)$ can also be carried out in the
spectral domain. Since  $\rho_u$ is
positive-definite, it constitutes, by Herglotz' theorem
\cite[p.~38]{KatzN},
the Fourier coefficients of a non-negative
measure $\mu := \mu_u$ on $\TT$ (the {\em power spectral measure}), i.e.,
\begin{equation}
\rho_u[k] = \int_\TT e^{-2\pi i k \xi} \,\dif \mu(\xi).
\end{equation}
Elementary Fourier analysis yields the spectral formula
\begin{equation} \label{MSE-spec}
\EE(x,\phi) = \int_\TT |2\sin(\pi\xi)|^{2m} |\Phi(\xi)|^2 \,\dif
\mu(\xi),
\end{equation}
where $\Phi$ has the absolutely convergent Fourier series representation
$$\Phi(\xi) = \sum_n \phi[n] e^{2\pi i n \xi}.$$
\par
This computational alternative is effective when the
measure $\mu$ has a simple description. On the other hand, it may happen
that this measure is somewhat complex to compute with directly.
As we shall demonstrate,
there will generally be a pure point (discrete) component
$\mu_\mathrm{pp}$ (i.e. a weighted sum of Dirac masses),
and an absolutely continuous component $\mu_\mathrm{ac}$
yielding a spectral density $s(\cdot)\in L^1(\TT)$,
where $\dif \mu_\mathrm{ac} (\xi) = s(\xi)\dif \xi$.
The continuous singular component will be nonexistent.
We shall analyze these two components
via their Fourier coefficients.
Under certain conditions, we will
be able to describe both components explicitly and
compute either asymptotics or sharp bounds for $\phi = \phi_M$
as $M\to \infty$.

\section{Piecewise affine maps of \sd quantization}
\label{piecewise-affine}

\par
In this section, we study the difference equations of \sd modulation
as a dynamical system arising from the iteration of
certain piecewise affine maps on $\RR^m$.
It easily follows from the first two equations in (\ref{state-eq-1}) that
\begin{equation}\label{ui[n]}
u_j[n]=\sum_{i=1}^j u_i[n{-}1] + (x[n]-q[n]), \;\;\;\;\;  1 \leq j \leq m,
\end{equation}
or in short,
\begin{equation}\label{Urec1}
{\bf u}[n]={\bf L}{\bf u}[n{-}1] + (x[n]-q[n]) {\bf 1},
\end{equation}
where ${\bf L}:={\bf L}_m$ is the $m{\times}m$ lower triangular matrix
of $1$'s and ${\bf 1} := {\bf 1}_m := (1,\dots,1)^\top \in \RR^m$.
Using the definition of
$q[n]$, we introduce a one-parameter family of maps
$\{\MM_x\}_{x\in \RR}$ on $\RR^m$ defined by
\begin{equation}\label{Mdef}
\textstyle
\MM_x ({\bf v}):={\bf L}{\bf v} + (x-Q(x,{\bf v})){\bf 1}.
\end{equation}
Hence, the evolution of the state vector ${\bf u}[n]$
is given by
\begin{equation}\label{recursive-u}
{\bf u}[n] =\MM_{x[n]} ({\bf u}[n{-}1]).
\end{equation}
According to the formulation presented in the introduction, the
elements of the partition $\Pi_x$ are then given by $\Omega_{x,i}
= \{ {\bf v}\in\RR^m ~:~ Q(x,{\bf v}) = d_i \}$, and the expression
\eqref{Mdef} is equivalent to \eqref{Mx}.
For the rest of the paper, we shall assume that $x[n]=x$ is a constant
sequence so that
\begin{equation}\label{recursive-u-constant-x}
{\bf u}[n] =\MM^{n}_{x} ({\bf u}[0]),
\end{equation}
and
\begin{equation}\label{recursive-q-constant-x}
q[n] =Q\!\left(x,\MM^{n-1}_{x} ({\bf u}[0])\right).
\end{equation}

\par
A variety of choices for the quantizer $Q$ have been
introduced in the practice of \sd modulation. Most of these
are designed with circuit implementation in mind, and therefore necessitate
simple arithmetic operations, such as linear combinations and
simple thresholding. A canonical
example would be
\begin{equation}\label{lin-rule}
Q_0(x,{\bf v}) = \lfloor \alpha_0x + \alpha_1v_1 + \dots + \alpha_mv_m +
\beta_0 \rfloor + \beta_1,
\end{equation}
where the coefficients $\alpha_i$ and $\beta_i$ are specific to each scheme.
We will call these rules ``linear'', referring to the fact that the
sets $\Omega_{x,i}$ are separated by translated hyperplanes in $\RR^m$.
There has also been recent research on
more general quantization rules and their benefits \cite{ingrid_devore,
SinanThao1, exp_decay}.
\par
Typically, an electrical
circuit cannot handle arbitrarily large amplitudes, and clips off
quantities that are beyond certain values. This is called
{\em overloading}. In this case, the effective mapping $Q$ is given by
\begin{equation}\label{overloaded}
Q(x,{\bf v}) = \left\{
\begin{array}{lll}
Q_0(x,{\bf v}) & \mbox{ if } & Q_0(x,{\bf v}) \in \{d_1,\dots,d_K\}, \\
d_1   & \mbox{ if } & Q_0(x,{\bf v}) < d_1, \\
d_K   & \mbox{ if } & Q_0(x,{\bf v}) > d_K.
\end{array}
\right.
\end{equation}
For the rest of the paper, 
we assume that the $d_i$ form a subset of an arithmetic progression of
spacing $1$, such as the case for the rule (\ref{lin-rule}). Since
we can always subtract a fixed constant from $x$ and the $d_i$,
we also assume, without loss of generality, that the $d_i$ are simply integers.
We shall be most interested in one-bit quantization rules,
i.e., rules for which $\mathrm{Ran}(Q) = \{d_1,d_2\}$.
Let us mention that one-bit \sd modulators are usually
overloaded by their nature.

\par
Let us emphasize once again that the quantization rule is crucial 
in the stability of the system. For a given $x$,
we call a \sd scheme defined by the quantization rule $Q(x,\cdot)$ {\em orbit
stable}, or simply {\em stable},
if for every initial condition ${\bf u}[0]$ in an
open set, the forward trajectory under the map $\MM_x$ is bounded in $\RR^m$,
and {\em positively stable}, if there exists a bounded 
set $\Gamma_0 \subset \RR^m$ with nonempty interior
that is positively invariant under $\MM_x$, i.e.,
$\MM_x(\Gamma_0) \subset \Gamma_0$. These two notions are
closely related. Clearly, positive stability
implies stability. On the other hand, in a stable scheme, 
if the forward trajectories of points in an open set 
are bounded with a uniform bound, then this 
would also imply the existence of a positively invariant bounded set. 
In practice, it is also desirable that stability holds uniformly in $x$.
However, we shall not need this kind of uniformity in this paper.

\par
\begin{figure}[t]
\begin{center}
\includegraphics[height=3 in]{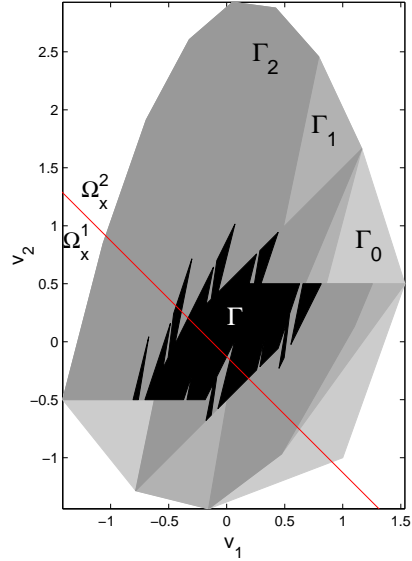}
\end{center}
\caption{\label{attractor}
The decreasing family of nested sets $\Gamma_k = \MM_x^k(\Gamma_0)$
indicated by decreasing
brightness. The limit set $\Gamma$ is invariant (see Theorem \ref{tiling1}).}
\end{figure}

In Figure \ref{attractor}, we depict a
positively invariant set $\Gamma_0$ under the map $\MM_x$ which
is defined by a one-bit linear rule in $\RR^2$. The set $\Gamma_0$
was found by a computer algorithm. In general, constructing positively
invariant sets for these maps is a non-trivial task \cite{Schreier,Yilmaz}.
Despite the presence of a vast collection of
\sd schemes that are used in hardware, only a small
set of them are proved to be stable. Most of the engineering
practice relies on extensive numerical simulation.

\par
In Figure \ref{attractor}, we also show in decreasing
brightness the forward iterates of $\Gamma_0$ given by
$\Gamma_k = \MM^k_x(\Gamma_0)$.
These sets converge to a limit set $\Gamma$, or the {\em attractor}, which
is shaded in black.
These invariant sets are the topic of discussion of next section.

\par
To avoid heavy and awkward notation, we shall drop the real parameter
$x$ from our notation except when we need it for a specific purpose
or for emphasis. It must be understood, however, that unless noted
otherwise, all objects that are derived from these dynamical systems
generally depend on $x$.

\section{Stability implies tiling invariant sets}
\label{sec-tiling}

\par
In this section we prove a crucial property
of the dynamics involved in positively stable \sd schemes. This is
called the {\em tiling property} and refers to the fact that there exist
trapping invariant sets that are disjoint unions of a
finite collection of disjoint tiles in $\RR^m$. Here a tile, or
a $\ZZ^m$-tile, means any subset $S$ of $\RR^m$ with the property that
$\{S + {\bf k} \}_{{\bf k} \in \ZZ^m}$ is a partition of $\RR^m$. Later in the
paper, this property will lead us to an exact spectral analysis
of the mean square error when the multiplicity of tiling is one.

\par
We consider a slightly more general class of
piecewise affine maps $\MM :=\MM_x$ on $\RR^m$, which are defined by
\begin{equation}
\MM({\bf v}) = \As_{x,i}({\bf v})
:= {\bf L}{\bf v} + x{\bf 1} + {\bf d}_i \;\;\;\;
\mbox{ if } \;\; {\bf v} \in \Omega_{x,i},
\end{equation}
where ${\bf L}$ is the lower triangular matrix of all $1$'s, and
$\{ \Omega_{x,i} \}_{i=1}^K$ is a finite Lebesgue measurable
partition of $\RR^m$, and ${\bf d}_i \in \ZZ^m$ for all $i=1,\dots,K$.
When ${\bf d}_i = -d_i {\bf 1}$, these maps are the same as those
that arise from \sd quantization.

\begin{theorem}\label{tiling1}\cite{ThaoBreaks}
Assume that there exists a bounded set $\Gamma_0 \subset \RR^m$ that
is positively invariant
under $\MM$, i.e., $\MM(\Gamma_0) \subset \Gamma_0$.
Then, the set $\Gamma \subset \Gamma_0$ defined by
\begin{equation}\label{gamma-def}
\Gamma := \bigcap_{k \geq 0} \MM^k(\Gamma_0)
\end{equation}
satisfies the following properties:
\begin{enumerate}
\item[(a)] $\MM(\Gamma) = \Gamma$,
\item[(b)] if $\Gamma_0$ contains a tile, then so does $\Gamma$.
\end{enumerate}
\end{theorem}

\begin{proof}%[Proof of Theorem \ref{tiling}]
This was previously proved in \cite{ThaoBreaks}.
For completeness of the discussion, we include the proof here.

\par {(a)}
Clearly, $\MM(\Gamma) \subset \Gamma \subset \Gamma_0$ since $\Gamma_0$
is positively invariant. We need to show that $\Gamma \subset \MM(\Gamma)$.
Let ${\bf v} \in \Gamma$ be an arbitrary point. Define
$\Gamma_k := \MM^k(\Gamma_0)$, $k\geq 0$. The sets $\Gamma_k$
form a decreasing sequence, and so is the case for the sets
$F_k := \MM^{-1}({\bf v}) \cap \Gamma_k$. Note that $\MM^{-1}({\bf v})$
is always finite since there are only finitely many
$\As_{x,i}$'s in the definition of $\MM$, each of which is 1-1.
($F_k$ would be finite even if there were infinitely many sets $\Omega_{x,i}$
because inverse images under $\MM$
have to differ by points in $\ZZ^m$ and only finitely many of them can
be present in $\Gamma_k$.) On the other hand ${\bf v} \in \Gamma_{k+1} =
\MM(\Gamma_k)$, therefore ${\bf v}$ has an
inverse image in $\Gamma_k$, i.e., $F_k$ is non-empty. Since $F_k$ form
a decreasing sequence of non-empty finite sets,
it follows that $\MM^{-1}({\bf v}) \cap \Gamma = \bigcap_{k\geq 0}
\,F_k \not= \emptyset$, i.e., ${\bf v} \in \MM(\Gamma)$. Hence
$\Gamma \subset \MM(\Gamma)$.

\par {(b)} Let $\Gamma_0$ contain a tile $G_0$, and define $G_k = \MM^k(G_0)$.
Each $G_k$ is a tile. To see this, note that for any given $i$,
$\As_{x,i}$ maps tiles to tiles, and for all
${\bf v} \in \RR^m$, $\MM({\bf v}) - \As_{x,i}({\bf v}) \in \ZZ^m$
so that $\MM$ maps tiles to tiles as well.
For an arbitrary point ${\bf w} \in \RR^m$, define the decreasing
sequence of sets $H_k = (\ZZ^m + {\bf w}) \cap \Gamma_k$.
Because $\Gamma_0$ is bounded, each $H_k$ is finite. On the other
hand, $\Gamma_k \supset G_k$ implies that each $\Gamma_k$ contains a tile,
yielding $H_k \not= \emptyset$.
Hence $(\ZZ^m + {\bf w}) \cap \Gamma =
\bigcap_{k\geq 0} H_k \not= \emptyset$. Since ${\bf w}$ is arbitrary, this
means that $\Gamma$ contains a tile.
\end{proof}
\par
In what follows, measurable means Lebesgue measurable, and
$\Leb(S)$ denotes the Lebesgue measure of a set $S$.

\begin{theorem}\label{tiling2}
Under the condition of Theorem \ref{tiling1}, assume moreover that
$x$ is irrational and that $\Gamma_0$ is measurable of non-zero
measure. Then, the set $\Gamma$ defined in \eqref{gamma-def}
differs from the union of a finite and non-empty collection of
disjoint $\ZZ^m$-tiles by a set of measure zero.
\end{theorem}

\begin{proof}
Clearly, $\Gamma$ is measurable since $\MM$ is piecewise affine.
Let us show that Lebesgue measure on $\Gamma$ is invariant under $\MM$.
From now on, we identify $\MM$ with its restriction on $\Gamma$.
From Theorem \ref{tiling1}, $\MM(\Gamma) = \Gamma$ which implies 
$\MM^{-1}(\Gamma) = \Gamma$ as well.
Let us define $A$ to be the set of points in $\Gamma$ with more than
one pre-image. A is measurable, simply because
$$ A = \bigcup_{i\not=j} \MM(\Gamma \cap \Omega_i) \cap \MM(\Gamma \cap
\Omega_j).$$
We claim that $\Leb(A) = 0$. Definition of $\MM$ implies
that $\MM$ preserves the measure of sets on which it is 1-1.
Since $\MM$ is 1-1 on $\MM^{-1}(\Gamma\backslash A)$, we have
$\Leb(\Gamma \backslash A) = \Leb(\MM^{-1}(\Gamma\backslash A))$.
On the other hand, since each point in $A$ has at least $2$ pre-images,
we have $2\Leb(A) \leq \Leb(\MM^{-1}(A))$. This implies
$$ 2\Leb(A) \leq \Leb(\MM^{-1}(A)) = \Leb(\MM^{-1}(\Gamma))
- \Leb(\MM^{-1}(\Gamma\backslash A)) = \Leb(\Gamma) -
\Leb(\Gamma \backslash A) = \Leb(A).$$
Therefore $\Leb(A) = \Leb(\MM^{-1}(A)) = 0$.
Hence, for any $B \subset \Gamma$, the
disjoint union $B = (B \cap A) \cup (B \backslash A)$ yields
$$\Leb(\MM^{-1}(B)) = \Leb(\MM^{-1}(B \cap A)) + \Leb(\MM^{-1}(B \backslash A))
= \Leb(B \backslash A) = \Leb(B),$$
i.e., $\MM$ preserves Lebesgue measure on $\Gamma$.

\par
Let $\pi:\Gamma \to \TT^m$ be the projection
defined by $\pi({\bf v}) = \langle {\bf v} \rangle$. Here we identify
$[0,1)^m$ with $\TT^m$. Let $\nu$ be the transformation of the measure
$\Leb|_\Gamma$ on $\TT^m$ under the projection
$\pi$, which is defined on the Lebesgue measurable subsets of $\TT^m$
by $\nu(B) = \Leb(\pi^{-1}(B))$. Let $\LL = \LL_x$ be
a generalized skew translation on $\TT^m$ defined by
\begin{equation}\label{def-L}
\LL{\bf v} := {\bf Lv} + x{\bf 1}\;\;\;\;\mbox{ (mod 1)}.
\end{equation}
Note that $\pi \MM = \LL \pi$.
Hence, for any measurable $B \subset \TT^m$, we have
$$ \nu(\LL^{-1}(B)) = \Leb(\pi^{-1}\LL^{-1}(B)) =
\Leb(\MM^{-1} \pi^{-1}(B)) = \Leb(\pi^{-1}(B)) = \nu(B),
$$
i.e., $\nu$ is invariant under $\LL$.

\par
At this point, we note that when $x$ is irrational,
$\LL$ is {\sl uniquely ergodic}, i.e., there is a
unique normalized non-trivial measure invariant under $\LL$,
which, in this case, is the Lebesgue measure. (See, for example,
\cite{Furst}, \cite[p.17]{Parry} for $m=2$, and
\cite[p.159]{Katok} for general $m$.\footnote{Here,
unique ergodicity is
stated for the map $(v_1,\dots,v_m) \mapsto
(v_1+x,v_2+v_1,\dots,v_m+v_{m-1})$, which is easily
shown to be isomorphic to $\LL$.})
Hence, $\vu = c\, \Leb$ for some $c \geq 0$; this includes the possibility
of the trivial invariant measure $\vu \equiv 0$.

\par
For each $j=0,1,\dots$, define
$$ T_j = \{ {\bf v} \in \TT^m ~:~ \mathrm{card}(\pi^{-1}({\bf v})) = j \}.$$
$\{T_j\}_{j \geq 0}$ is a finite measurable partition of $\TT^m$.
The finiteness is due to the fact that
$\Gamma$ is a bounded set and measurability is
simply due to the relation
$$ T_j = \left \{{\bf v} \in \RR^m ~:~ \sum_{{\bf k} \in \ZZ^m}
\chi_{_{\Gamma + {\bf k}}}({\bf v}) = j \right \}. $$
Note that
$$c\, \Leb(T_j) = \nu(T_j) = \Leb(\pi^{-1}(T_j)) = j \Leb(T_j).$$
This shows that there cannot exist two such sets $T_i$ and $T_j$ both
with non-zero measure. Hence, there exists a (unique) $j$, namely, $j=c$,
such that
$\Leb(\TT^m \backslash T_j) = 0$. This implies that $\Gamma$ is the union of
$j$ copies of $\TT^m$, possibly with the exception of a set of zero measure.

\par
Let us now show that $j\geq 1$.
Consider $\Sigma_0:=\pi(\Gamma_0) \subset \TT^m$. Since $\Gamma_0$ is
positively invariant, we find that $\LL(\Sigma_0) =
\pi \MM(\Gamma_0) \subset \pi(\Gamma_0) = \Sigma_0$. Since
$\LL$ is 1-1, we have $\LL^{-1}(\Sigma_0) \supset \Sigma_0$. Hence,
$ \LL^{-1}(\Sigma_0)\; \triangle\; \Sigma_0 =  
\LL^{-1}(\Sigma_0) \backslash \Sigma_0
= \LL^{-1}(\Sigma_0 \backslash \LL(\Sigma_0))$. This implies, since
$\LL$ is measure-preserving,
$$ \Leb(\LL^{-1}(\Sigma_0)\; \triangle \;\Sigma_0)
= \Leb(\LL^{-1}(\Sigma_0 \backslash \LL(\Sigma_0))) 
= \Leb(\Sigma_0 \backslash \LL(\Sigma_0)) 
= \Leb(\Sigma_0) - \Leb(\LL(\Sigma_0)) = 0.$$
Ergodicity of
$\LL$ implies that $\Leb(\Sigma_0)$ is 0 or 1. The first case is not
possible, since each point in $\Sigma_0$ has at most finitely many
inverse images under $\pi^{-1}$ and this would violate
$\Leb(\Gamma_0) > 0$. Therefore $\Leb(\Sigma_0) = 1$, implying that $j\geq 1$.
\end{proof}

\par When $x$ is irrational,
Theorem \ref{tiling2} improves Theorem \ref{tiling1}(b)
in two aspects. First, the outcome is that $\Gamma$ not only contains
a tile, but in fact it is {\sl composed} of disjoint tiles, up 
to a set of measure zero. Second, to conclude this,
it suffices to check that $\Gamma_0$ has positive measure, instead of
the stronger (though equivalent) requirement that $\Gamma_0$ contain a tile. 
On the other hand, Theorem \ref{tiling1}(b) is still 
interesting due to its algebraic nature: It can can be used to 
test if $\Gamma$ contains an exact tile (i.e., $\pi(\Gamma) = \TT^m$), and
it remains valid even when $x$ is rational.
\par
Let us also note, as an application of Theorem \ref{tiling2}, that
whenever a positively invariant set $\Gamma_0$ of $\MM_x$ (for irrational $x$)
can be found with $0 < \Leb(\Gamma_0) < 2$, then the invariant set $\Gamma$ 
is a single tile. 
\par
\begin{figure}[t]
\begin{center}
\includegraphics[height=3in]{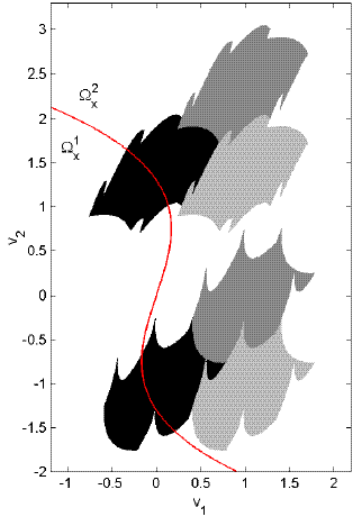}
\end{center}
\caption{\label{2-tiles}
Represented in black is the invariant set $\Gamma$ 
of a $1$-bit $2$nd order scheme whose partition
is determined by the cubic curve shown in the figure.
The copies in gray are the translated versions of $\Gamma$
by $(1,0)$ and $(1,1)$, respectively. In this example, each connected
component of $\Gamma$ is also invariant.
}
\end{figure}

 
In Figure \ref{2-tiles}, we show an illustration of an invariant
set which is composed of two tiles. In this example, the \sd scheme
is 1-bit 2nd order and the partition is determined by a cubic curve.

\section{The single-tile case and its consequence}

\label{conseq-tiling} 
\begin{figure}[t]
\centerline{\hbox{
\vbox{
 \hbox{\makebox[6.0cm]{\includegraphics[height=5cm]{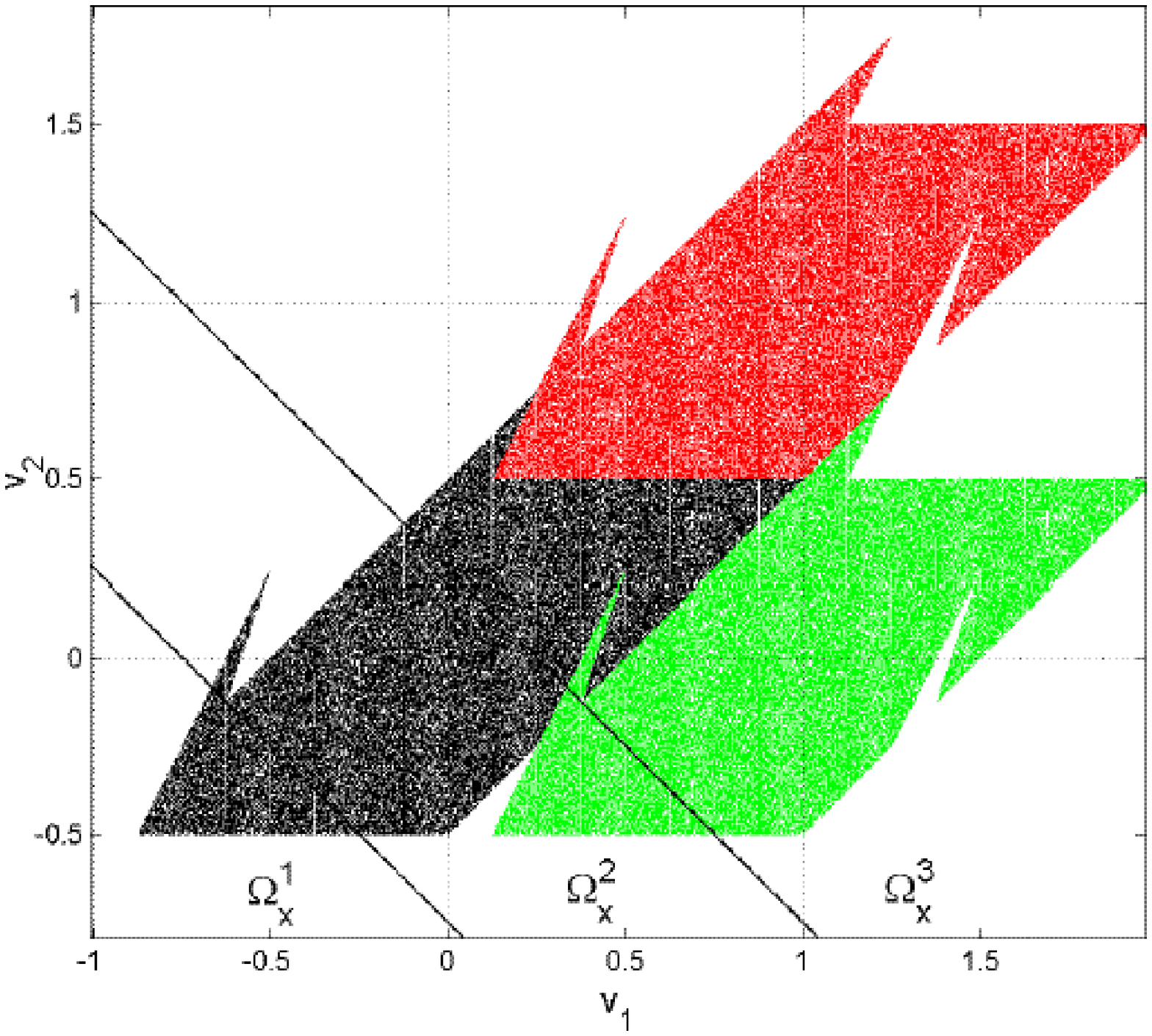}
}}
 \hbox{\makebox[6.0cm]{(a)}} }
\vbox{\hbox{\makebox[1mm]{}}} \vbox{
\hbox{\makebox[6.0cm]{\includegraphics[height=5cm]{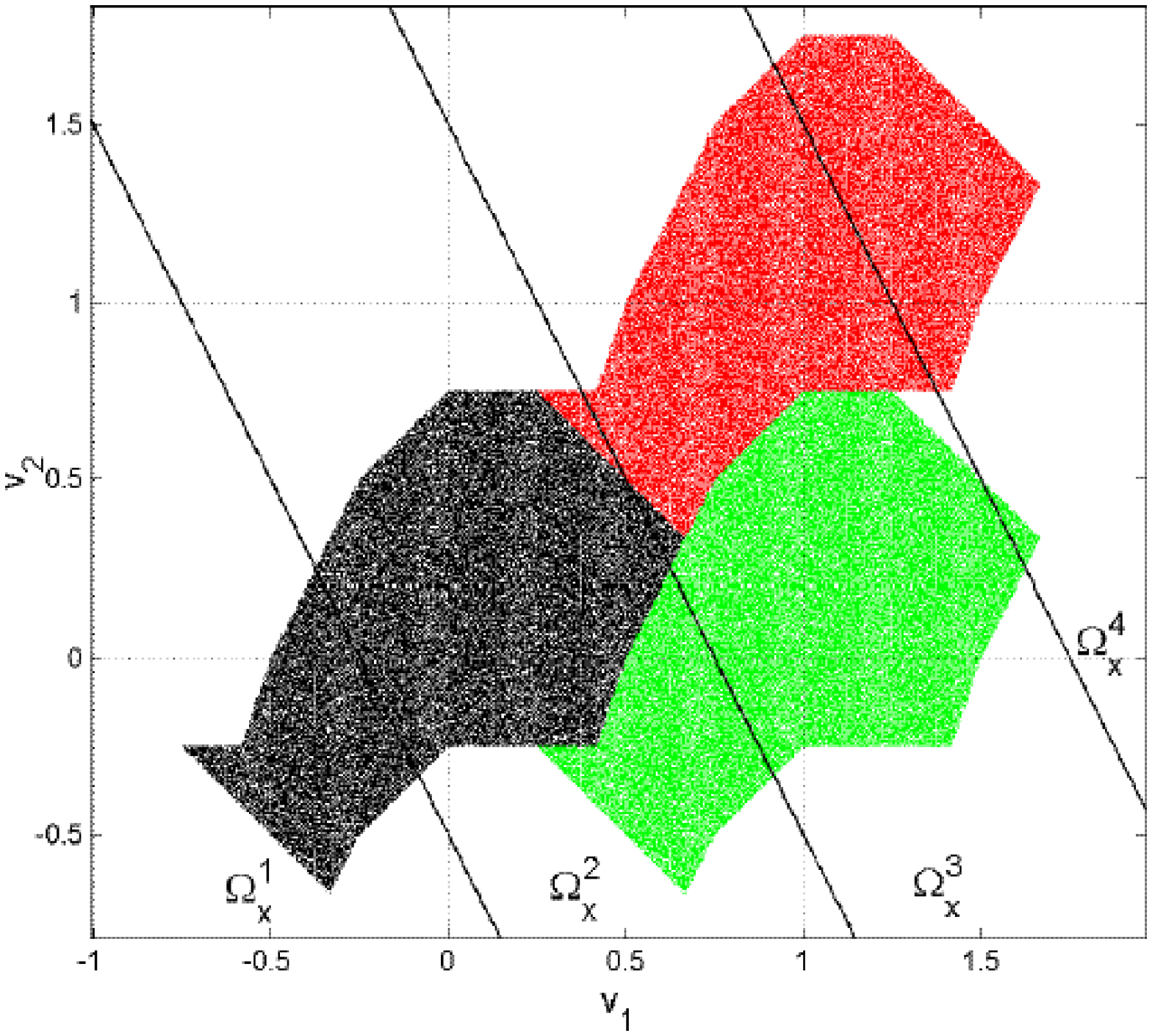}
 }}
 \hbox{\makebox[6.0cm]{(b)}} }
  }}
\centerline{\hbox{
\vbox{
 \hbox{\makebox[6.0cm]{\includegraphics[height=5cm]{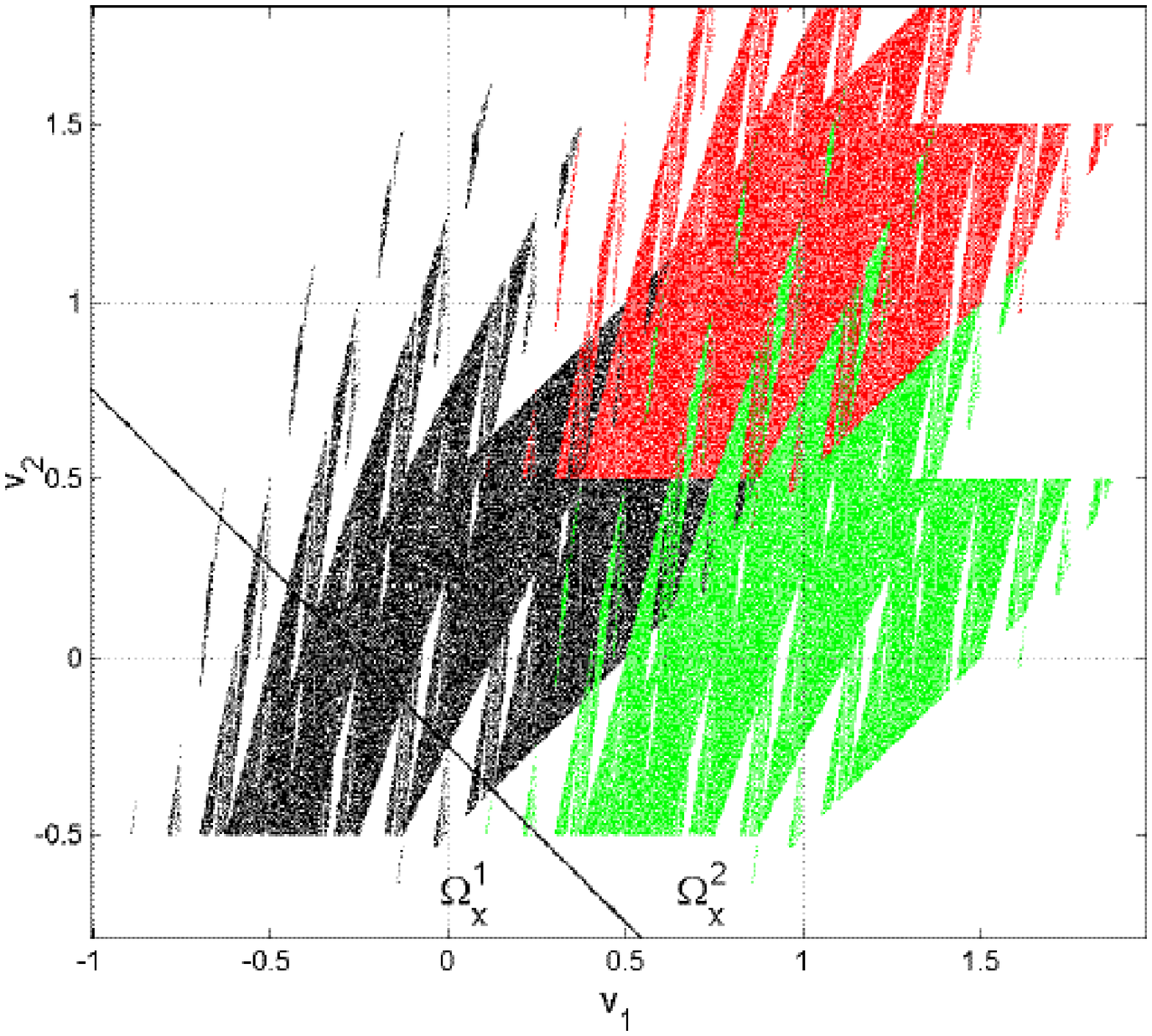}
}}
 \hbox{\makebox[6.0cm]{(c)}} }
\vbox{\hbox{\makebox[1mm]{}}} \vbox{
\hbox{\makebox[6.0cm]{\includegraphics[height=5cm]{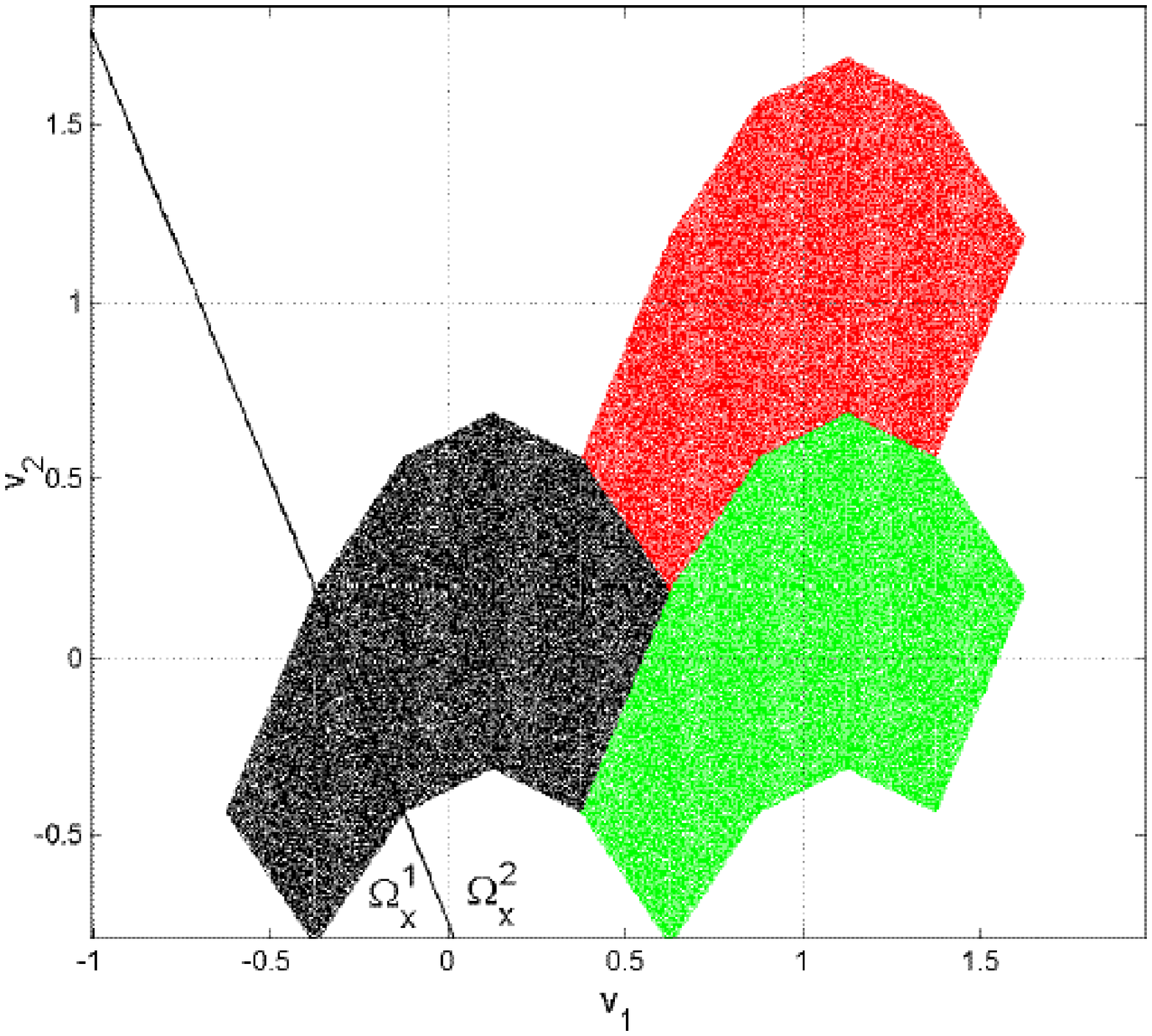}
}}
 \hbox{\makebox[6.0cm]{(d)}} }
  }}
\caption{\label{exp-tiling}
Representation in black of several consecutive state points ${\bf u}[n]$
of various second
order $\Sigma\Delta$ modulators with the irrational input $x \approx 3/4$.
The copies in gray are the translated versions of the state points by $(1,0)$
and $(1,1)$, respectively.}
\end{figure}

\begin{figure}[htbp]
\centerline{\hskip 0.45cm \includegraphics[height=2in]{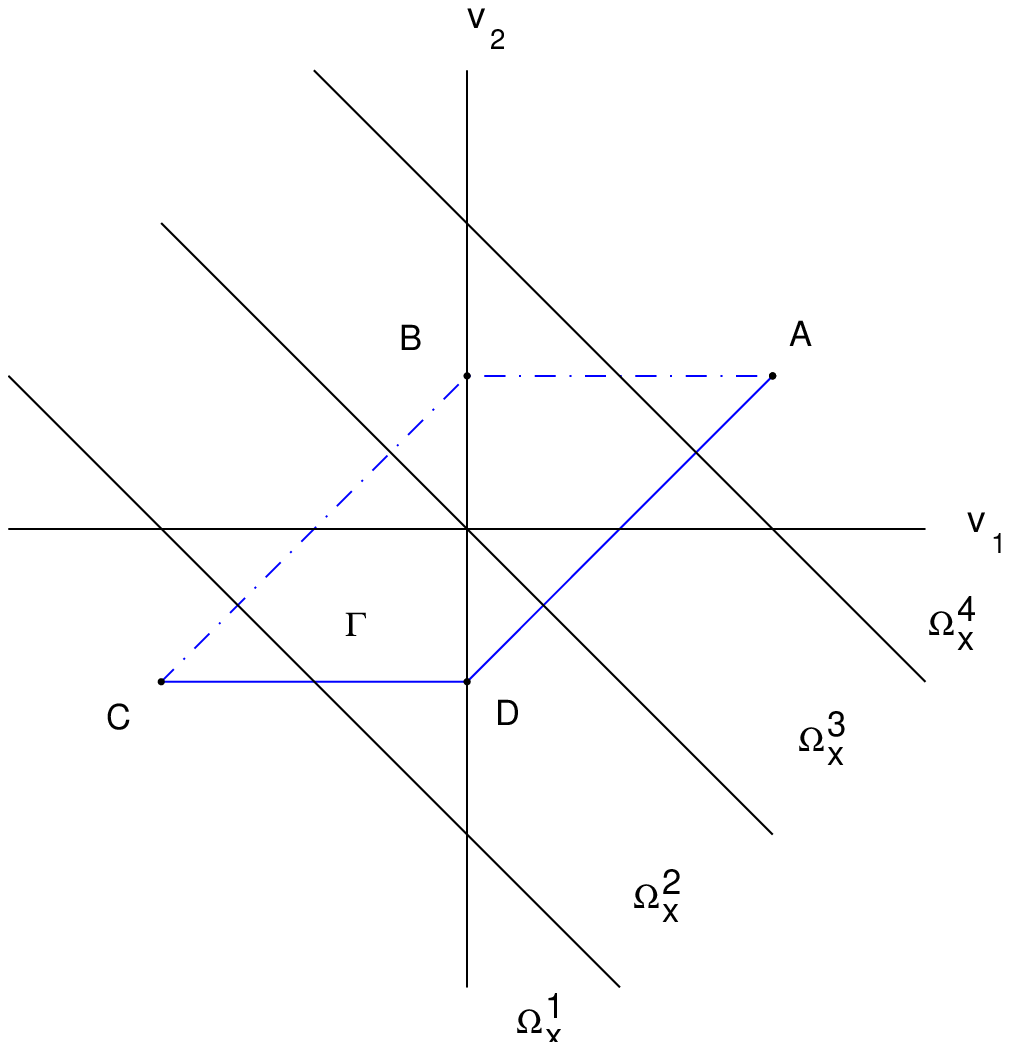}}
\centerline{(a) 2-bit ``linear'' with $(d_1,d_2,d_3,d_4) =
(-1,0,1,2)$. ($x=0.5$)}

\vskip 1 cm

\centerline{
\includegraphics[height=1.9in]{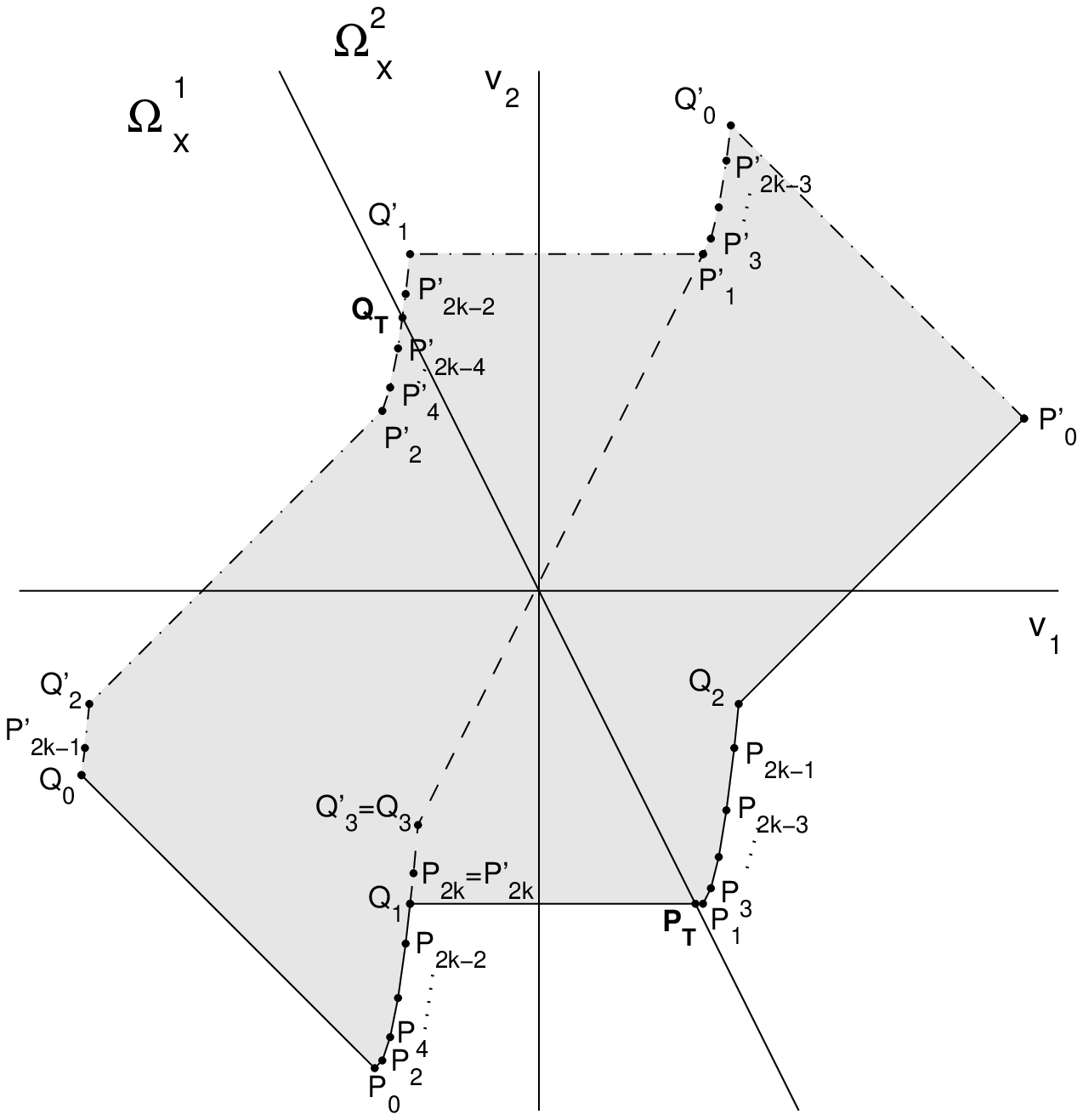} \hskip 1 cm
\includegraphics[height=1.9in]{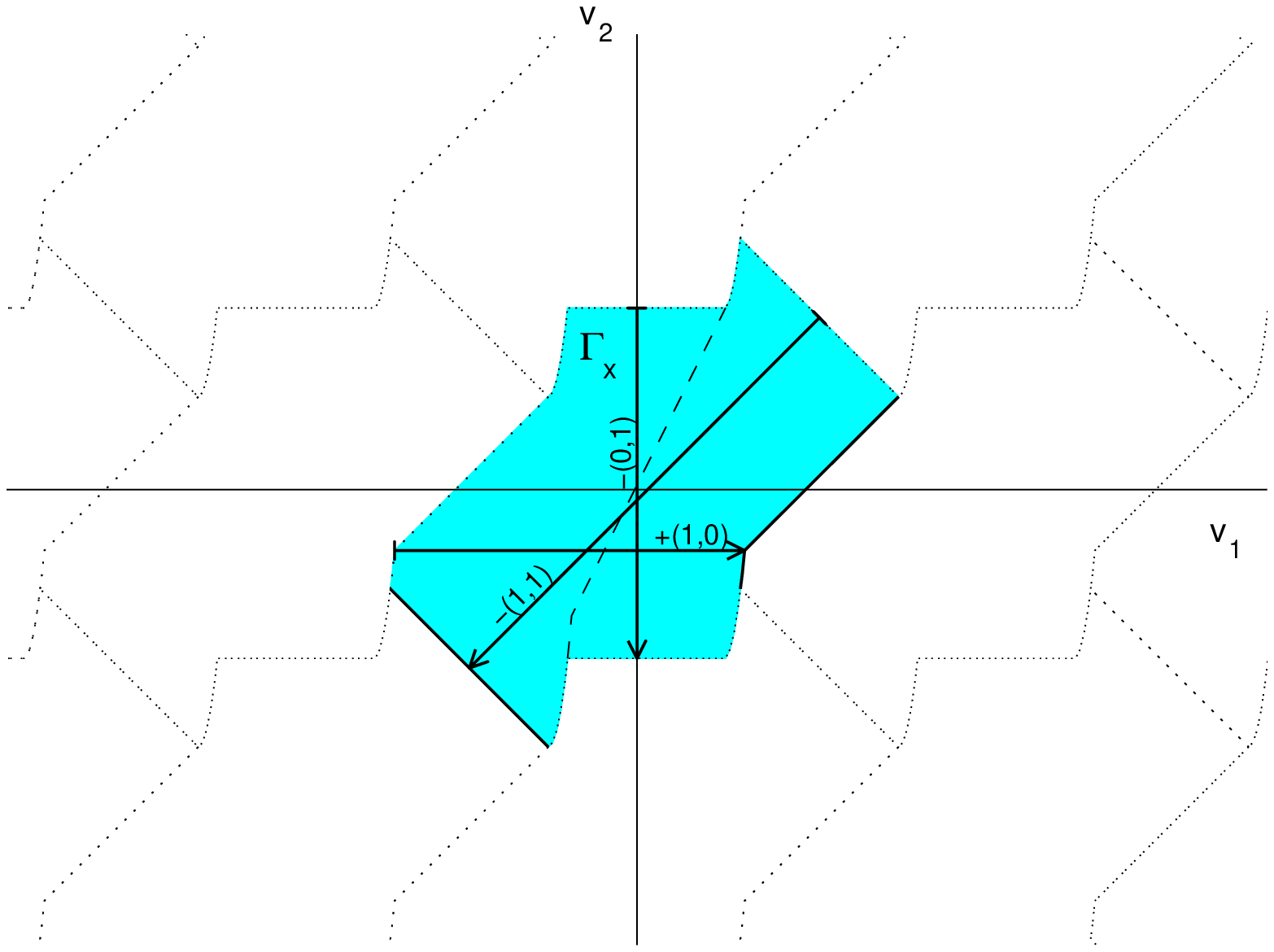}
}
\centerline{(b) 1-bit ``linear'' with $(d_1,d_2) =
(0,1)$. ($x \approx 0.52$)}

\vskip 1 cm

\centerline{\includegraphics[height=1.7in]{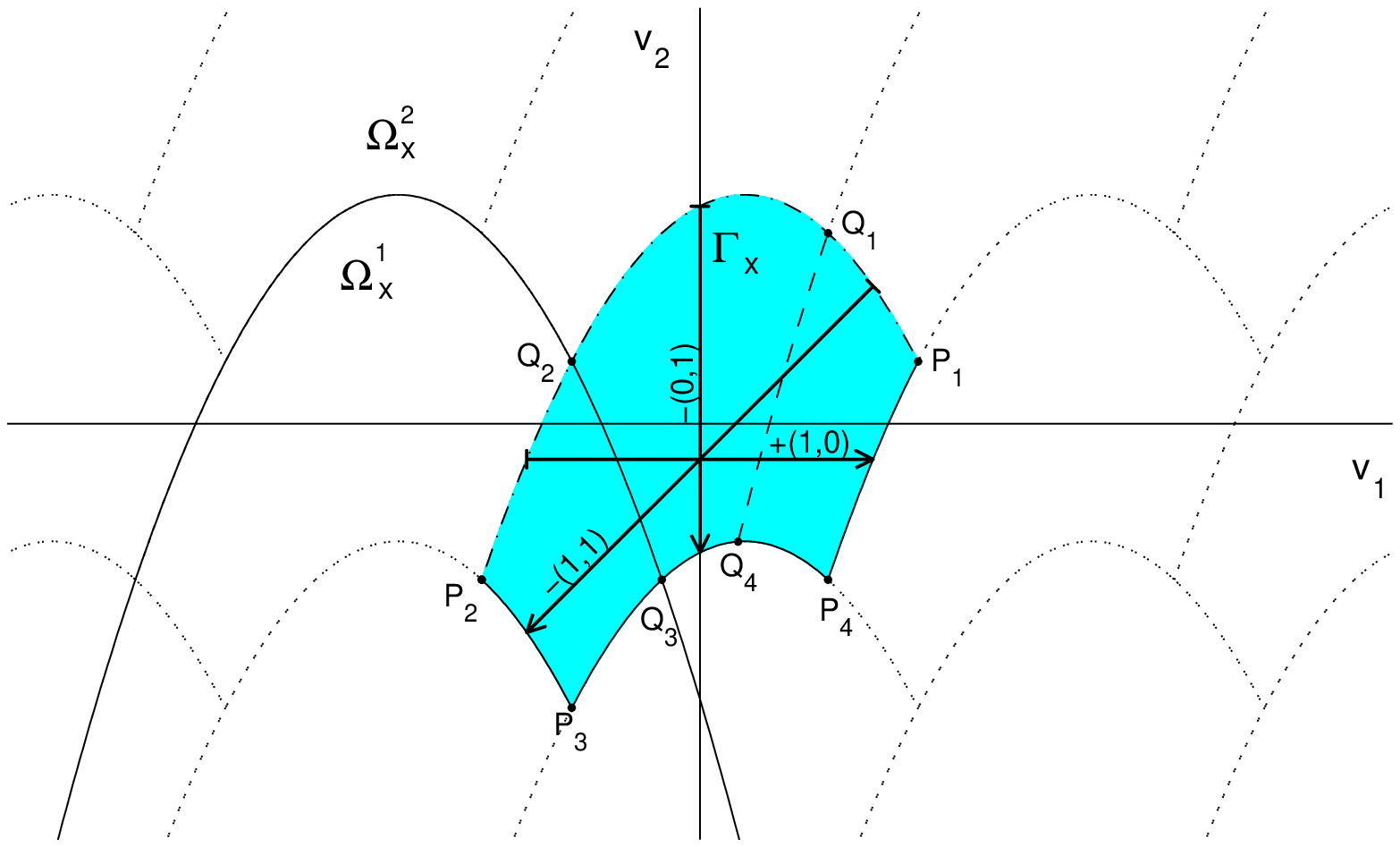}}
\centerline{(c) 1-bit ``quadratic'' with $(d_1,d_2) =
(0, 1)$. ($x = 0.74$)}
\caption{\label{proven-tiling}
Three families of quantization rules for which the tiling property was 
proven in \cite{SinanThao1} with parametric
explicit expressions for the corresponding invariant sets.}
\end{figure}

Since the initial experimental discovery of
the tiling property in \cite{CSG_thesis,SinanThao1}, we have observed that
the invariant
sets $\Gamma$ resulting from  practical stable second order \sd
schemes systematically appear to be single
tiles. We show in Figure \ref{exp-tiling} experimental examples of
$\Gamma$ on second order schemes. In Figure \ref{proven-tiling},
we show the set $\Gamma$ in three cases where its explicit
analytical derivation has been possible \cite{SinanThao1}. (In
these particular cases, $\Gamma$ is actually proven to be an exact
tile.) A fundamental question is to characterize maps $\MM_x$ which
yield a single invariant tile. In this paper, we will simply
assume that this condition is realized. As will be seen, a whole
new framework of analysis will be generated from this particular
situation.

\par
A tile $\Gamma$ intrinsically generates a unique projection
$\langle\cdot\rangle_{_{\Gamma}}: \RR^m \to \Gamma$
such that ${\bf v} {-} \langle {\bf v}\rangle_{_{\Gamma}}
\in \ZZ^m$ for all ${\bf v} \in \RR^m$.
The restriction of this $\ZZ^m$-periodic projection
to the unit cube $[0,1)^m$ (which
we identify with  $\TT^m$)
is a measure preserving bijection (note that the inverse of 
$\langle\cdot\rangle_{_{\Gamma}}:\TT^m \to \Gamma$
is the map $\pi$ that was defined in the proof of Theorem \ref{tiling2}). 
When $\Gamma$ is invariant
under $\MM$, the map $\langle\cdot\rangle_{_{\Gamma}}: \TT^m \to \Gamma$
establishes an isomorphism
between $\MM$ on $\Gamma$ and the affine transformation $\LL := \LL_x$
on $\TT^m$ defined by (\ref{def-L}).
Indeed, the definition of $\LL$ easily yields
$\LL ({\bf v}) {-} \MM
(\langle {\bf v}\rangle_{_{\Gamma}}) \in \ZZ^m$. Hence,
$$\langle \LL ({\bf v}) \rangle_{_{\Gamma}}
= \MM \!
\left(\langle {\bf v}\rangle_{_{\Gamma}}\right), $$
or in other words, the following diagram commutes:
\[
      \begin{CD}
      {\TT^m} @> {\LL} >> {\TT^m} \\
      @ V{\langle\cdot\rangle_{_{\Gamma}}} VV @
VV {\langle\cdot\rangle_{_{\Gamma}}} V \\
      {\Gamma} @>> {\MM} > {\Gamma} \\
      \end{CD}
     \]

%\begin{figure*}
%\centerline{\hbox{ \vbox{
%\hbox{\makebox[5mm]{(a)}} \vspace{2.3cm}} \vbox{
% \hbox{\makebox[6.6cm]{
%\psfig{figure=s_tiling2a.eps,height=5.6cm} }}}
%\vbox{\hbox{\makebox[5mm]{}}} \vbox{ \hbox{\makebox[5mm]{(b)}}
%\vspace{2.3cm}} \vbox{ \hbox{\makebox[6.6cm]{
%\psfig{figure=s_tiling2b.eps,height=5.6cm} }}} }}
%\centerline{\hbox{ \vbox{ \hbox{\makebox[5mm]{(c)}}
%\vspace{2.3cm}} \vbox{
% \hbox{\makebox[6.6cm]{
%\psfig{figure=s_tiling2c.eps,height=5.6cm} }}}
%\vbox{\hbox{\makebox[5mm]{}}} \vbox{ \hbox{\makebox[5mm]{(d)}}
%\vspace{2.3cm}} \vbox{ \hbox{\makebox[6.6cm]{
%\psfig{figure=s_tiling2d.eps,height=5.6cm} }}} }}
%\centerline{\hbox{ \vbox{ \hbox{\makebox[5mm]{(e)}}
%\vspace{2.3cm}} \vbox{
% \hbox{\makebox[6.6cm]{
%\psfig{figure=s_tiling2e.eps,height=5.6cm} }}}
%\vbox{\hbox{\makebox[5mm]{}}} \vbox{ \hbox{\makebox[5mm]{(f)}}
%\vspace{2.3cm}} \vbox{ \hbox{\makebox[6.6cm]{
%\psfig{figure=s_tiling2f.eps,height=5.6cm} }}} }}
%\caption{\label{tiling} Representation in black of 10,000
%consecutive state points ${\bf u}[k]$ of various second order \sd
%modulators with constant inputs $x[k]$: (a,c,e)
%$x[k]=x_0=\sqrt{2}/11$; (b,d,f) $x[k]=x_1=\sqrt{2}/4$; (a,b)
%multi-bit configuration with $A(z)=1-0.3z^{-1}$; (c,d) multi-bit
%configuration with $A(z)=1$ and a deviation of the 0-threshold by
%$1/3$; (e,f) single-bit configuration with $A(z)=1$. The points in
%gray are the translated versions of the state vector points by
%$[1~0]$ and $[1~1]$, respectively.}
%\end{figure*}

\par
The first important consequence of single invariant tiles is that it
reduces the dynamical system $\MM$ to
the much simpler $\LL$ whose $n$-fold composition can be
computed explicitly. It follows that
if ${\bf u}[0] \in \Gamma$, then
\begin{equation}\label{u-n-formula}
{\bf u}[n] = \MM^{n}({\bf u}[0]) =
\langle \LL^{n} ({\bf u}[0]) \rangle_{_{\Gamma}}
= \Big \langle {\bf L}^n {\bf u}[0] + x {\bf s}[n] \Big
\rangle_{_{\Gamma}},
\end{equation}
where ${\bf s}[n] := {\bf s}_m[n]$ is defined by
\begin{equation}
{\bf s}[n] = \left(\sum_{k=0}^{n-1} {\bf L}^{k} \right){\bf 1}.
\end{equation}
It is an easy computation to show that
the $j$th coordinate of ${\bf s}[n]$, which we denote by $s_j[n]$,
is equal to $\binom{j{+}n{-}1}{j}$.

\par
The second important consequence is that 
if $x$ is an irrational number, then $\MM_x$ on $\Gamma$
inherits the ergodicity of $\LL_x$ via the isomorphism generated
by $\langle\cdot\rangle_{_{\Gamma}}$. Since $\langle\cdot\rangle_{_{\Gamma}}
:\TT^m \to \Gamma$ preserves Lebesgue measure,  $\MM_x$ is then ergodic
with respect to the restriction of Lebesgue measure on $\Gamma$.
Hence Birkhoff Ergodic Theorem yields

\begin{proposition}\label{ergodic}
Let $x$ be an irrational number and
$\Gamma$ be a Lebesgue measurable $\ZZ^m$-tile 
(up to a set of measure zero) that is invariant under
$\MM$. Then for any function $F \in L^1(\Gamma)$,
\begin{equation}\label{unif-n-F}
\lim_{N\rightarrow\infty}\frac{1}{N}\sum_{n=1}^{N}F({\bf u}[n])=
\int_{\Gamma}F({\bf v})\,\dif {\bf v}
= \int_{\TT^m} F(\langle {\bf v}\rangle_{_{\Gamma}})\,\dif {\bf v}
\end{equation}
for almost every initial condition ${\bf u}[0] \in \Gamma$.
\end{proposition}
This formula will be the fundamental computational tool for the
analysis of the autocorrelation sequence $\rho_u$. 
For the remainder of this paper, we shall assume that we are working
with quantization rules for which the invariant sets are composed of
single tiles.
This will save us from repetition in the assumptions of our results.
However, it will also be important to know certain
geometric features of these invariant tiles.
We will state these explicitly when needed.

\section{Analysis of the autocorrelation sequence $\rho_{u}$}

\par
Let $\pr({\bf v}) = v_m$ be the projection of
a vector ${\bf v} \in \RR^m$ onto its $m$th coordinate. If we define
the function
\begin{equation}
F_k({\bf v}) = \pr({\bf v})\pr(\MM^k({\bf v})),
\end{equation}
then it follows that
$$u_m[n]u_m[n+k] = \pr({\bf u}[n])\pr(\MM^k({\bf u}[n])) =
F_k({\bf u}[n]),$$
and therefore Proposition \ref{ergodic}
gives an expression for the value of  $\rho_u[k]$:
\begin{equation}\label{auto-integral}
\rho_u[k] = \int_{\Gamma} F_k({\bf v})\,\dif {\bf v}
= \int_{\TT^m} F_k(\langle {\bf v}\rangle_{_{\Gamma}})\,\dif {\bf v}.
\end{equation}
A direct evaluation of $\rho_u[k]$
in either of these forms is not easy, because the
$k$-fold iterated map $\MM^k$
as well as the invariant set $\Gamma$ are implicitly-defined
and complex objects.
The problem can be somewhat simplified via the conjugate map
$\LL^k$. Indeed, one has
$$F_k \circ \langle \cdot \rangle_{_{\Gamma}}
= \Big(\pr \circ \langle \cdot \rangle_{_{\Gamma}}\Big)
\Big (\pr \circ \MM^k \circ \langle \cdot \rangle_{_{\Gamma}}
\Big )
=  \Big(\pr \circ \langle \cdot \rangle_{_{\Gamma}}\Big)
\Big (\pr \circ \langle \cdot \rangle_{_{\Gamma}} \circ \LL^k
\Big ),
$$
so that if we define
$$G_\Gamma = \pr \circ \langle \cdot \rangle_{_{\Gamma}},$$
then via (\ref{u-n-formula}), we obtain the formula
\begin{equation}\label{rho-G-formula-1}
\rho_u[k] = \int_{\TT^m} G_\Gamma({\bf v})G_\Gamma(\LL^k ({\bf v}))
\,\dif {\bf v} =
\int_{\TT^m} G_\Gamma({\bf v})G_\Gamma({\bf L}^k{\bf v}
+ x{\bf s}[k])\,\dif {\bf v},
\end{equation}
which now only depends on $\Gamma$.
\par
As it is standard in the spectral theory of dynamical systems
(see, e.g., \cite{Parry}),
let $\UU:=\UU_\LL$ be the unitary operator on $L^2(\TT^m)$
defined by $(\UU f)({\bf v}) = f (\LL({\bf v}))$.
Then \eqref{rho-G-formula-1} reduces to
\begin{equation}\label{rho-G-formula-2}
\rho_u[k] = \left(G_\Gamma,\, \UU^k G_\Gamma\right)_{L^2(\TT^m)}.
\end{equation}
For any $f \in  L^2(\TT^m)$, the inner products
$\left(f,\, \UU^k f\right)_{L^2(\TT^m)}$, $k\in \ZZ$, define a
positive-definite sequence
so that there exists a unique non-negative measure
$\nu_f$ on $\TT$ with Fourier coefficients
\begin{equation}
\hat \nu_f[k] = \left(f,\,\UU^k f\right)_{L^2(\TT^m)}
\end{equation}
for all $k \in \ZZ$.
Note that when $f = G_\Gamma$, it follows from (\ref{rho-G-formula-2}) that
the corresponding measure $\nu_{G_\Gamma} = \mu$, where $\mu$ is
the spectral measure that
was mentioned in Section \ref{sec-error}, with $\hat \mu = \rho_u$.

\subsection{Decomposition of the mixed spectrum: General results}

\par
We shall separate the autocorrelation sequence $\rho_u$ into two
additive components that result from two different types of
spectral behavior. Using the spectral theorem for unitary
operators, we decompose $ L^2(\TT^m) $ into two $\UU$-invariant,
orthogonal subspaces as $L^2(\TT^m) = \HH_\mathrm{pp} \oplus
\HH_{\mathrm c}$, where
$$\HH_\mathrm{pp} = \{f \in L^2(\TT^m) : \nu_f \mbox{ is
purely atomic} \},$$ which is also equal to
the closed linear span of the set of
all eigenfunctions of $\UU$, and
$$ \HH_{\mathrm c} = \HH_\mathrm{pp}^\perp =
\{f \in L^2(\TT^m) : \nu_f \mbox{ is non-atomic (continuous)} \}.$$
\par
In the particular case of the transformation $\LL$,
it turns out that every spectrum on $\HH_\mathrm{pp}^\perp$ is
absolutely continuous (see Appendix A).
Therefore we denote $ \HH_\mathrm{c}$ by $ \HH_\mathrm{ac}$.
Any $f \in L^2(\TT^m)$ can be uniquely decomposed as
$f = f_\mathrm{pp} + f_\mathrm{ac}$,
where $f_\mathrm{pp} \in \HH_\mathrm{pp}$ and
$f_\mathrm{ac} \in \HH_\mathrm{ac}$.
For $\LL$, it is also known, as we show in Appendix A, that
$$
\HH_\mathrm{pp} =
\{f \in L^2(\TT^m) : f({\bf v}) \mbox{ only depends on } v_1 \},
$$
and the orthogonal projection of $f$ onto $\HH_\mathrm{pp}$ is
given by
\begin{equation}
f_\mathrm{pp}({\bf v}) =
\int_{\TT^{m-1}} f(v_1,{\bf v}') \,\dif {\bf v}'.
\end{equation}
For notational convenience, we
write $\bar f :=f_\mathrm{pp}$ and $\breve f := f_\mathrm{ac}$
for any $f \in L^2(\TT)$.
\par
For  $f = G_\Gamma$, we now consider the decomposition
\begin{equation}
G_\Gamma = \bar{G}_\Gamma + \breve{G}_\Gamma.
\end{equation}
Because of orthogonality and $\UU$-invariance
of $\HH_\mathrm{pp}$ and $\HH_\mathrm{ac}$ , (\ref{rho-G-formula-2}) implies that
\begin{equation}
\rho_u[k] = \left(\bar G_\Gamma, \,\UU^k \bar G_\Gamma\right)_{L^2(\TT^m)} +
\left(\breve G_\Gamma, \,\UU^k \breve G_\Gamma\right)_{L^2(\TT^m)},
\end{equation}
providing the decomposition
$$\rho_u =  \bar \rho_u + \breve \rho_u.$$
Here, using formula
(\ref{u-n-formula}) and the fact that functions in the subspace
$\HH_\mathrm{pp}$ depend only on the first variable, we obtain
\begin{equation}
\bar \rho_u[k] =
\left(\bar G_\Gamma, \,\UU^k \bar G_\Gamma\right)_{L^2(\TT^m)}
= \int_{\TT} \bar G_\Gamma(v_1) \bar G_\Gamma(v_1 + k x) \,\dif v_1
\end{equation}
and
\begin{equation}
\breve \rho_u[k] =
\left(\breve G_\Gamma, \,\UU^k \breve G_\Gamma\right)_{L^2(\TT^m)}
= \int_{\TT^m} \breve G_\Gamma({\bf v}) \breve G_\Gamma
({\bf L}^k{\bf v} + x {\bf s}[k]) \,\dif {\bf v}.
\end{equation}
This decomposition provides the Fourier coefficients of the
pure-point $\mu_\mathrm{pp}$ and the absolutely
continuous $\mu_\mathrm{ac}$ components of the spectral measure, respectively.
It also yields an explicit simple formula for
$\mu_\mathrm{pp}$ in terms of the
Fourier coefficients of $\bar G_\Gamma$. We have

\begin{theorem} \label{pure-point}
\begin{equation}\label{mu-pp}
\mu_\mathrm{pp} = \sum_{n \in \ZZ}
\left| \widehat{\bar G_\Gamma}[n] \right|^2 \,
\delta_{nx},
\end{equation}
where $\delta_a$ denotes the unit Dirac mass at $a \in \TT$.
\end{theorem}

\begin{proof}
Let $\nu$ denote the measure given on the right hand side
of (\ref{mu-pp}).
It suffices to verify that $\hat \nu[k]=\bar \rho_u[k]$
for all $k\in \ZZ$. We find by direct evaluation that
$$ \hat \nu[k] = \sum_{n \in \ZZ}
\left| \widehat{\bar G_\Gamma}[n] \right|^2 e^{-2\pi i k nx}.
$$
Clearly, this is the Fourier expansion of the even-symmetric function
$$ A_{\bar G_\Gamma}(\xi) = \int_\TT \bar G_\Gamma(v) \bar G_\Gamma(v+\xi)
\, \dif v$$
evaluated at $\xi = -kx$. Therefore $\hat \nu[k] =
\bar \rho_u[-k] = \bar \rho_u[k]$.
\end{proof}

\par \noindent 
{\bf Note:} It is easy to see that this result holds for any function
$f \in L^2(\TT^m)$ in the sense that
\begin{equation}
(\nu_f)_\mathrm{pp} = \sum_{n \in \ZZ}
\Big| \widehat{f}_\mathrm{pp}[n] \Big|^2 \, \delta_{nx}.
\end{equation}

\par On the other hand, the computation of
$\mu_\mathrm{ac}$ is not easy. Since absolute continuity implies an
integrable density $s(\cdot)$, where
$d\mu_\mathrm{ac} (\xi) = s(\xi)\dif \xi$,
we can immediately say by the Riemann-Lebesgue lemma that
the Fourier coefficients $\breve \rho_u[k] \to 0$ as $|k| \to \infty$.
However, the rate of decay is determined by the geometry of $\Gamma$. We
shall be particularly interested in the case when the density
is continuous at $\xi = 0$.

\subsection{Properties of $\breve \rho_u$ for
the class of $v_m$-connected invariant tiles}

\par In this section, we derive explicit formulae for the Fourier
coefficients $\breve \rho_u[k]$ of the continuous
component of the spectral measure $\mu$ when the
invariant tile $\Gamma$ has certain geometric regularity.
For a given  tile $\Gamma$ for $\RR^m$, let us define
\begin{equation}
\Lambda_\Gamma := \bigcup_{\quad{\bf k}' \in \ZZ^{m-1}}
\Gamma + ({\bf k}',0),
\end{equation}
and for any ${\bf v}' \in \RR^{m-1}$,
\begin{equation}
\Lambda_\Gamma({\bf v}') := \pr(\Lambda_\Gamma \, \cap \,
\{{\bf v}'\}{\times}\RR) =
\Big  \{ v_m \in \RR : ({\bf v}', v_m) \in \Lambda_\Gamma \Big \}.
\end{equation}

\begin{proposition}\label{lambda-set}
For each ${\bf v}'\in \RR^{m-1}$, the set
$\Lambda_\Gamma({\bf v}')$ is a tile
in $\RR$ with respect to $\ZZ$-translations, and
\begin{equation}
G_\Gamma({\bf v}',v_m) = \langle v_m \rangle_{\Lambda_\Gamma({\bf v}')}.
\end{equation}
\end{proposition}
\begin{proof}
Since $\Gamma$ is a tile, the collection of sets
$\{ \Lambda_\Gamma + ({\bf 0}, k) : k \in \ZZ \}$ forms a partition
of $\RR^m$. Therefore for any  ${\bf v}'\in \RR^{m-1}$,
the $v_m$-section of this collection given by
$\{\Lambda_\Gamma({\bf v}') + k : k \in \ZZ \}$,
is a partition of $\RR$. This shows that
$\Lambda_\Gamma({\bf v}')$ is a tile.
For the second part of the claim, let ${\bf v}=({\bf v}',v_m)$.
The definition of $\pr$ immediately yields
$$({\bf v}',G_\Gamma({\bf v})) = ({\bf v}',\pr(\langle {\bf v}
\rangle_{_\Gamma}))
= \langle {\bf v}\rangle_{_\Gamma} + ({\bf k}',0)$$
for some ${\bf k}' \in \ZZ^{m-1}$. This says that
$({\bf v}',G_\Gamma({\bf v})) \in \Lambda_\Gamma$ and
therefore $G_\Gamma({\bf v}) \in \Lambda_\Gamma({\bf v}')$.
The result follows since $G_\Gamma({\bf v}',v_m){-}v_m  \in \ZZ$.
\end{proof}
\begin{definition} We say that a tile $\Gamma \subset \RR^m$ is
$v_m$-connected if for each ${\bf v}'\in \RR^{m-1}$, the one dimensional
tile $\Lambda_\Gamma({\bf v}')$ is a connected set, i.e. a unit-length
interval. In this case, we denote by $\lambda_\Gamma({\bf v}')$
the midpoint of $\Lambda_\Gamma({\bf v}')$.
\end{definition}

\par
\begin{figure}[htbp]
\centerline{\hbox{
\vbox{
\hbox{\makebox[5mm]{(a)}}
\vspace{2.3cm}}
\vbox{
 \hbox{\makebox[7.5cm]{
\includegraphics[height=4.5cm]{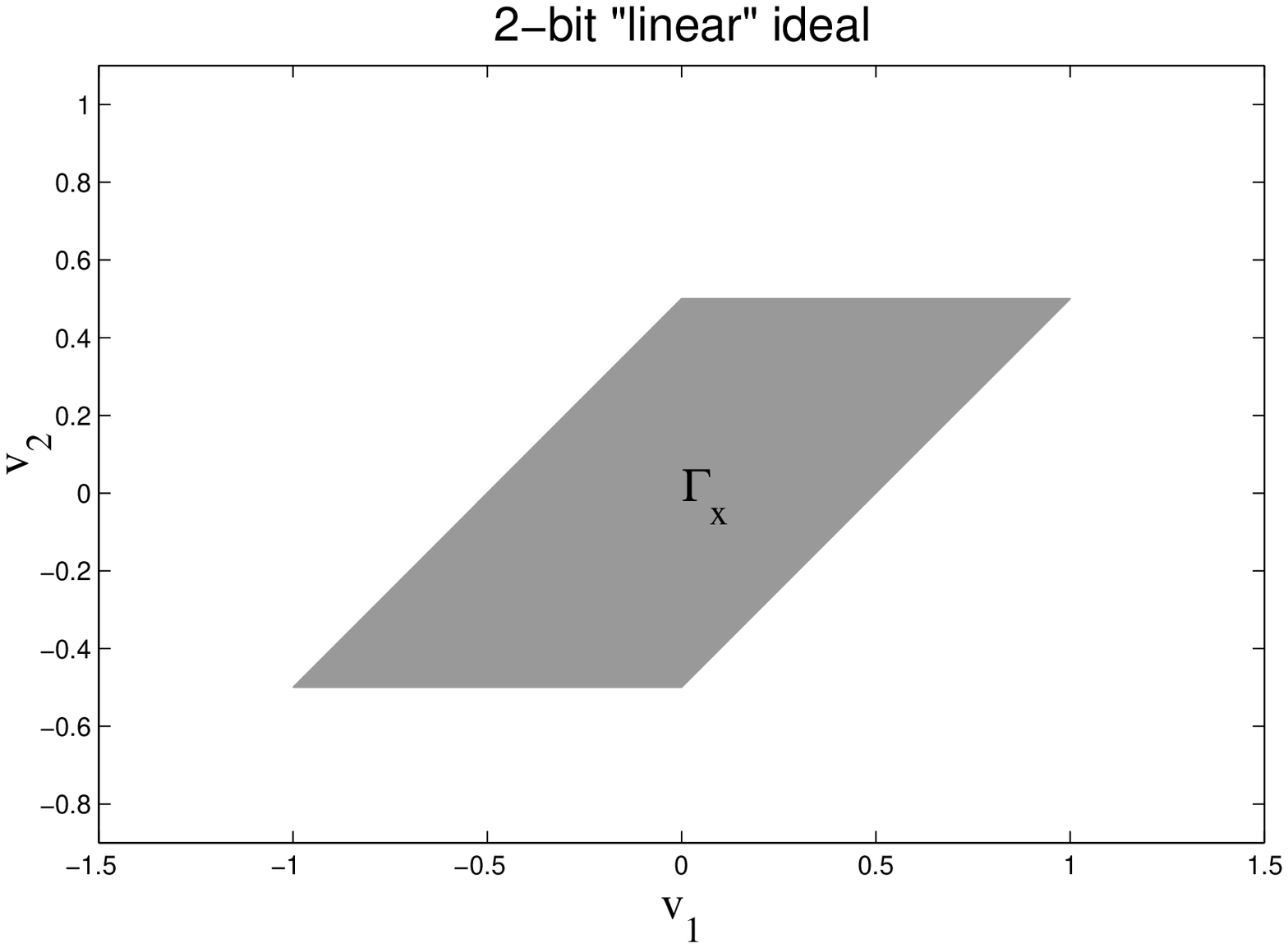}
}}} \vbox{
\hbox{\makebox[7.5cm]{ 
\includegraphics[height=4.5cm]{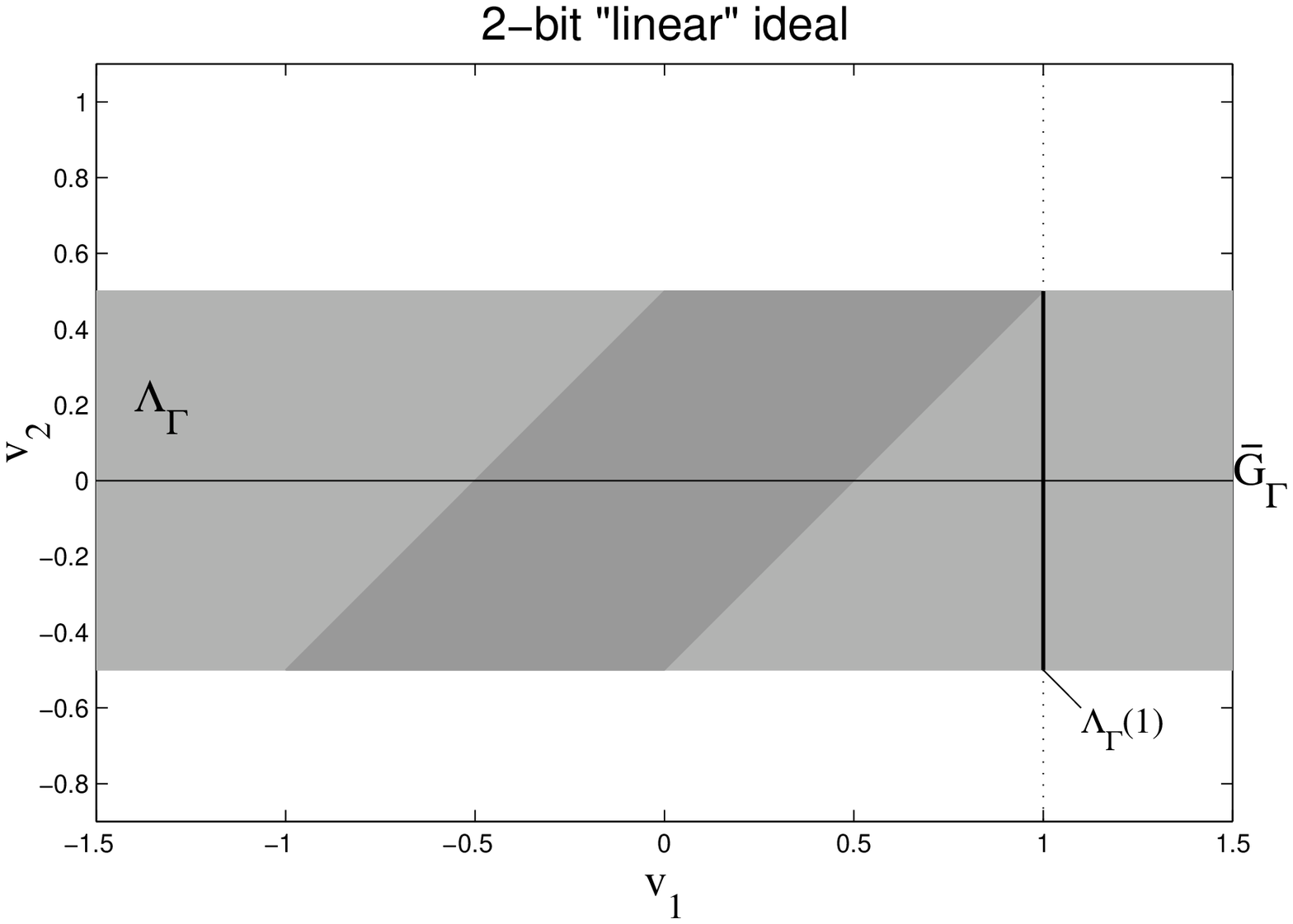}
}}}
}} \centerline{\hbox{ \vbox{ \hbox{\makebox[5mm]{(b)}}
\vspace{2.3cm}} \vbox{
 \hbox{\makebox[7.5cm]{
\includegraphics[height=4.5cm]{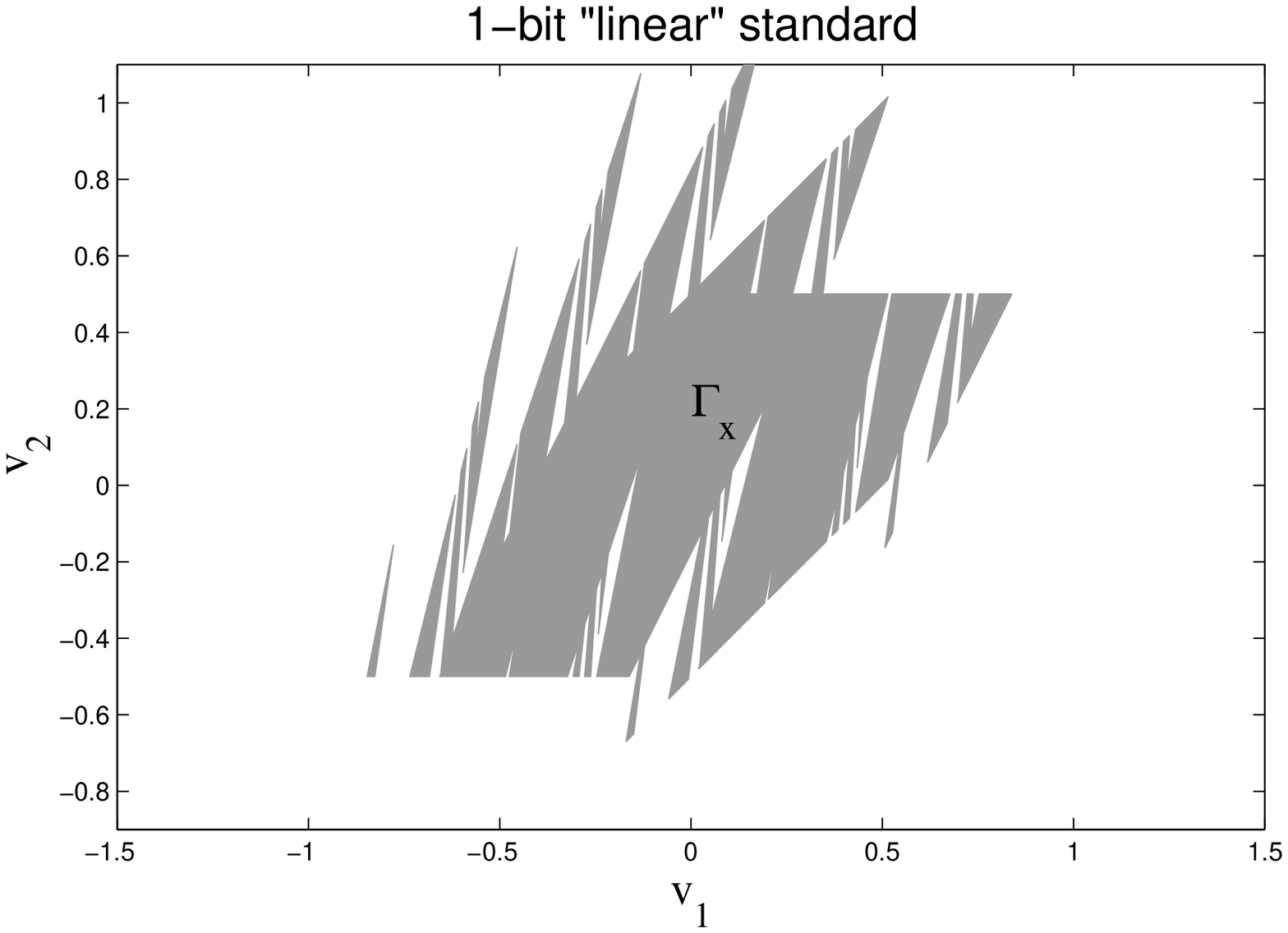}
}}} \vbox{
\hbox{\makebox[7.5cm]{ 
\includegraphics[height=4.5cm]{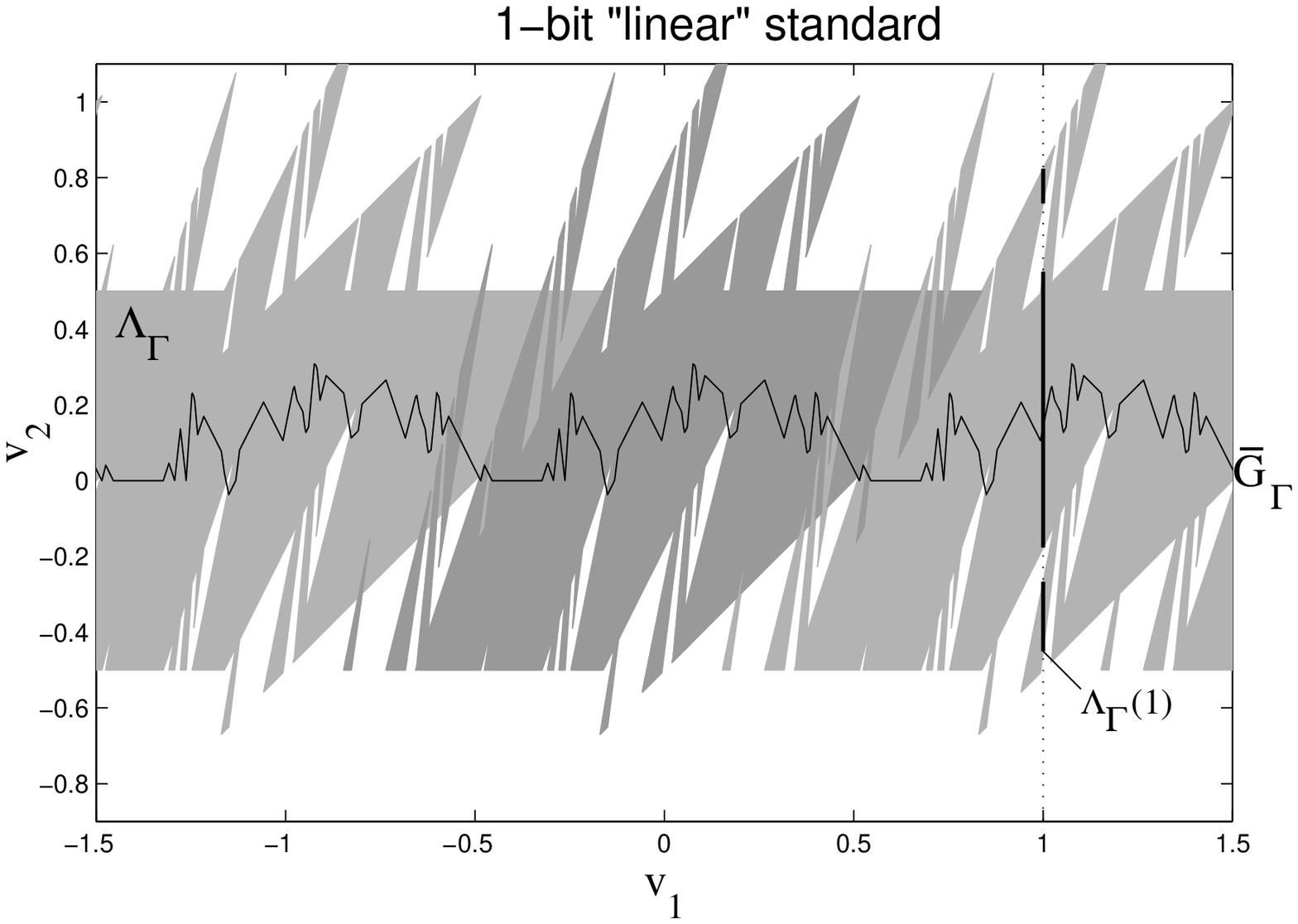}
}}}
}} \centerline{\hbox{ \vbox{ \hbox{\makebox[5mm]{(c)}}
\vspace{2.3cm}} \vbox{
 \hbox{\makebox[7.5cm]{
\includegraphics[height=4.5cm]{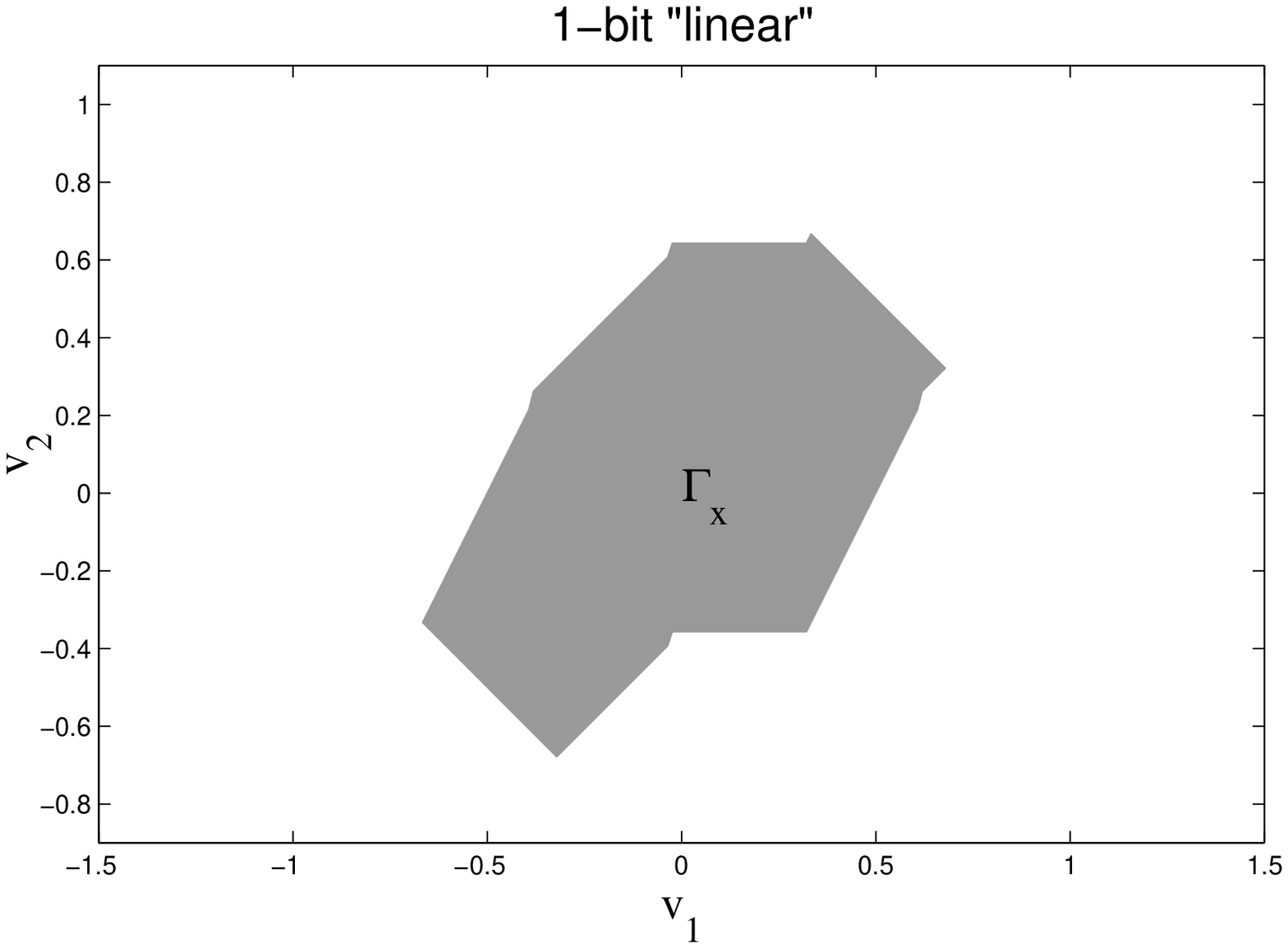}
}}} \vbox{
\hbox{\makebox[7.5cm]{ 
\includegraphics[height=4.5cm]{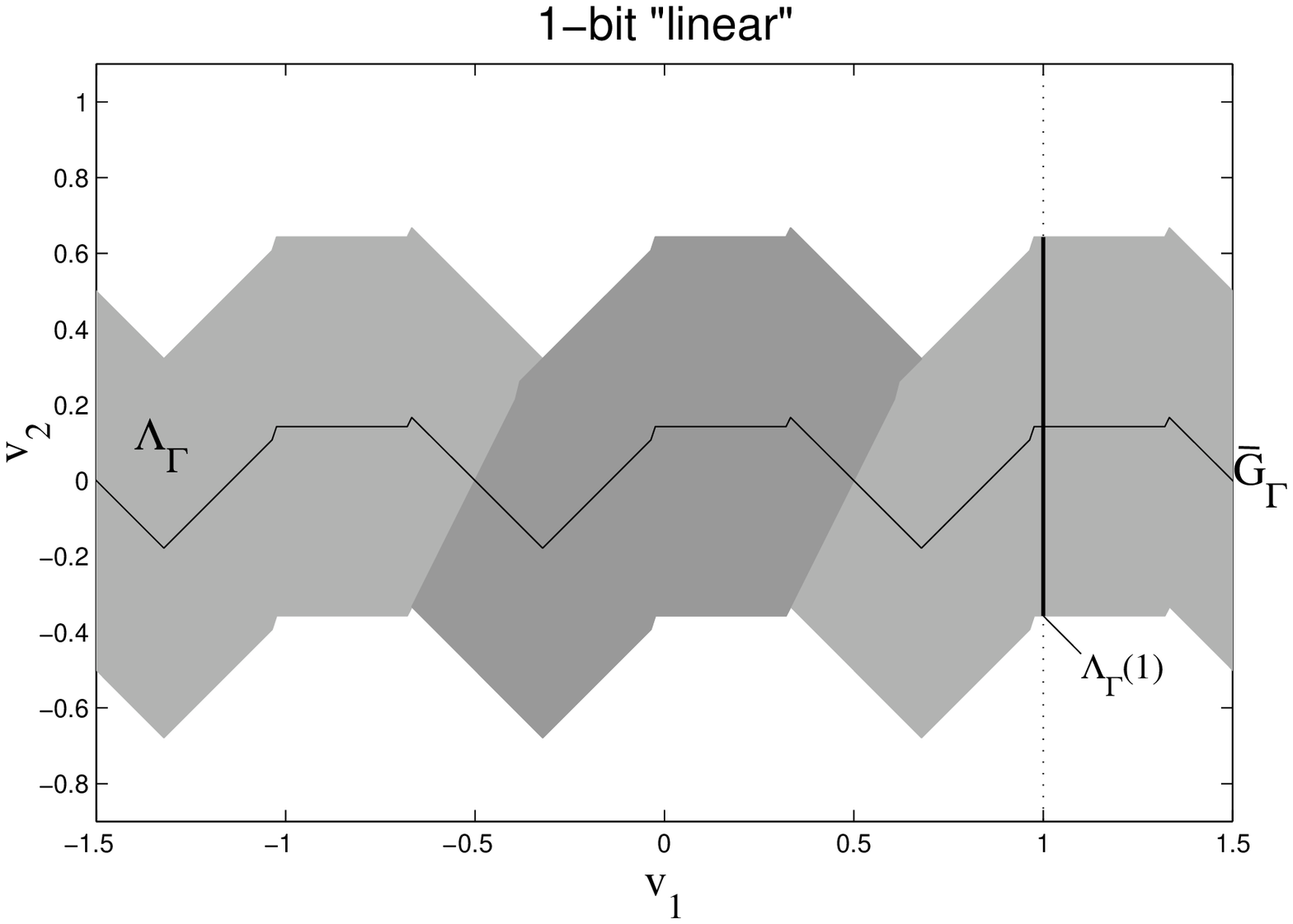}
}}}
}} \centerline{\hbox{ \vbox{ \hbox{\makebox[5mm]{(d)}}
\vspace{2.3cm}} \vbox{
 \hbox{\makebox[7.5cm]{
\includegraphics[height=4.5cm]{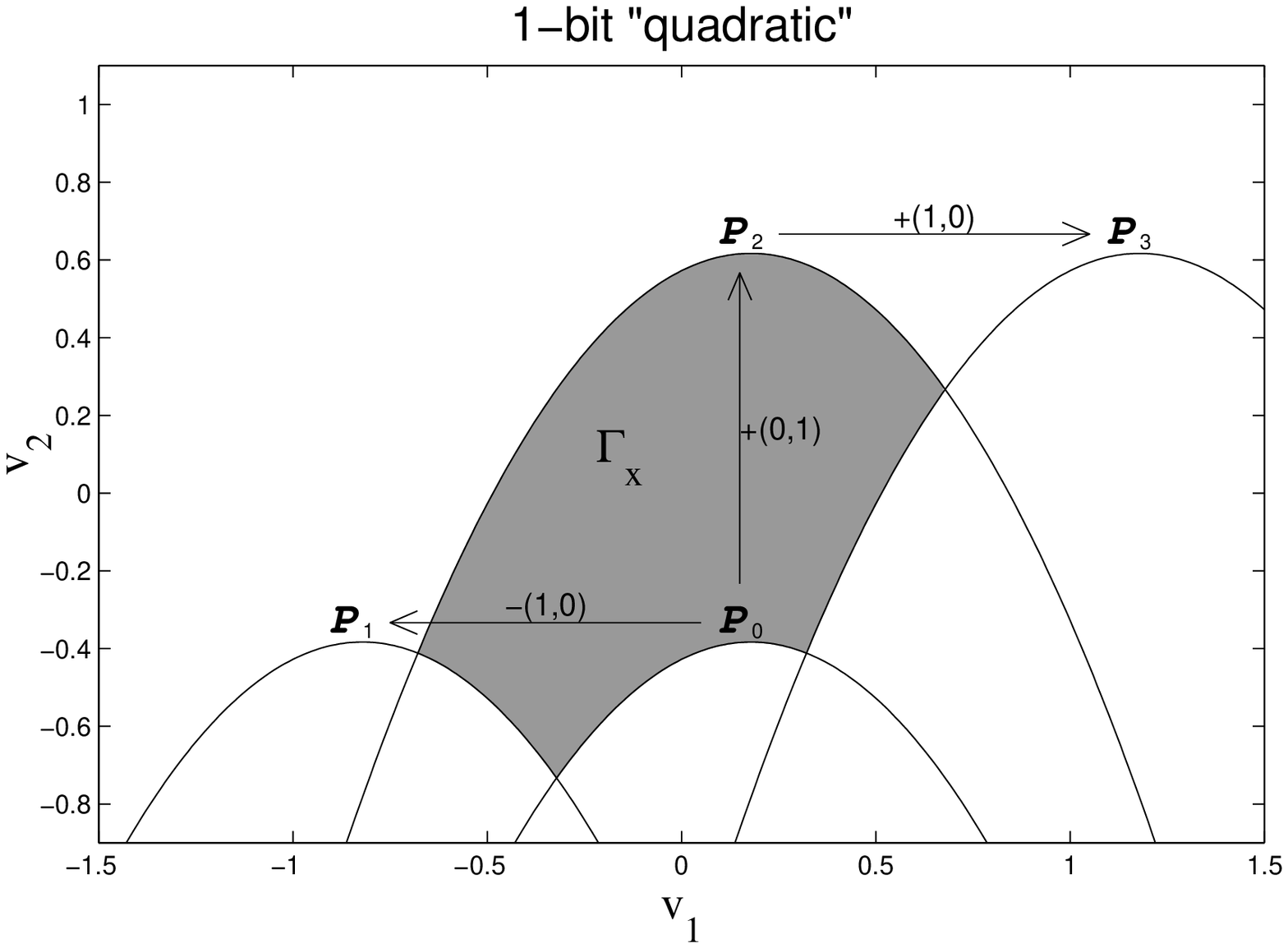}
}}} \vbox{
\hbox{\makebox[7.5cm]{ 
\includegraphics[height=4.5cm]{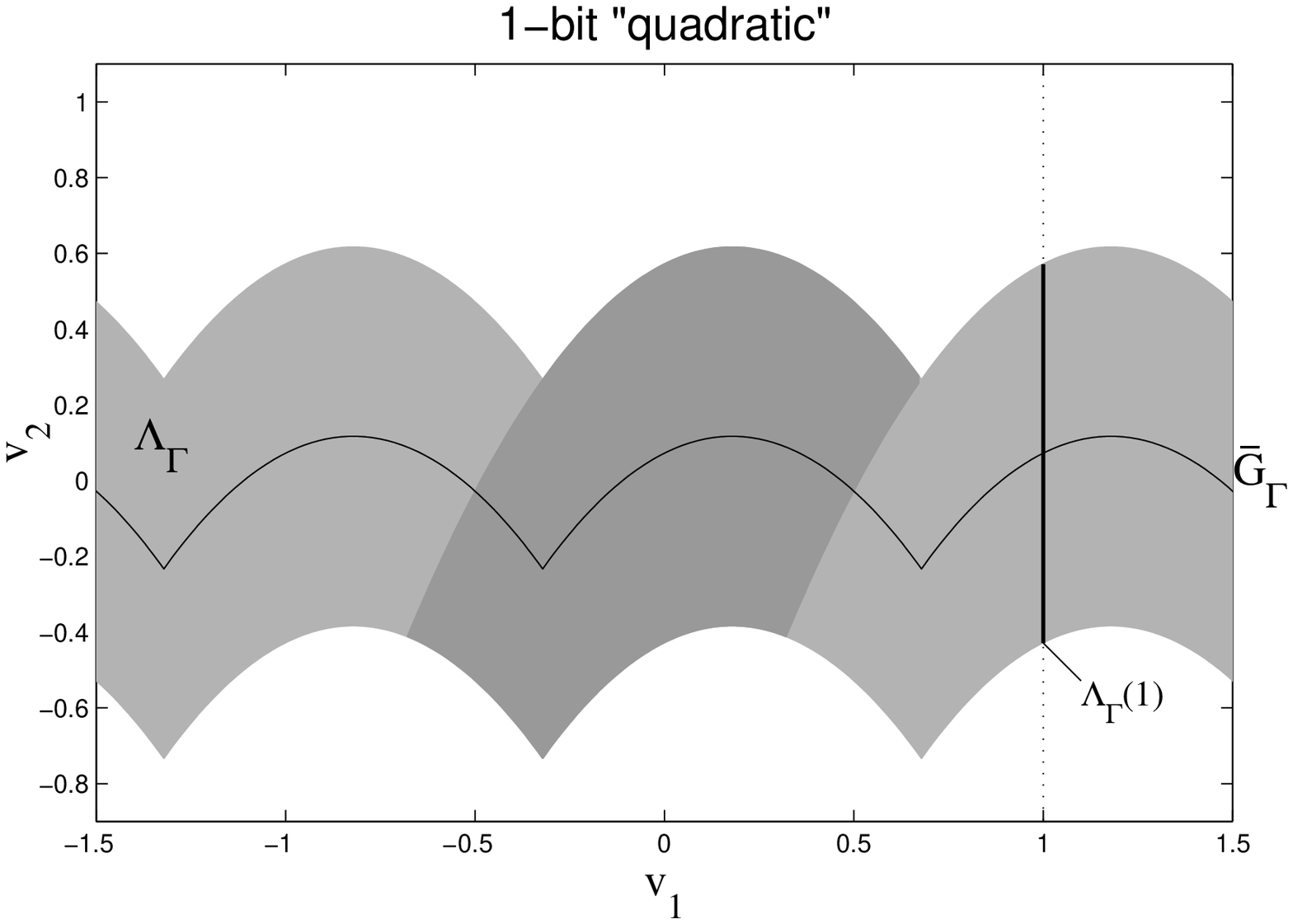}
}}}
}} \centerline{\hbox{ \vbox{ \hbox{\makebox[5mm]{}}}
\vbox{\hbox{\makebox[7.5cm]{(i)}}} \vbox{
\hbox{\makebox[7.5cm]{(ii)}}} }} \caption{\label{invar} Invariant
tiles of various second order modulators: (i) Invariant tile
$\Gamma_x$, (ii) corresponding set $\Lambda_\Gamma$.}
\end{figure}

 In Figure \ref{invar},
we display examples of the function $\Lambda_\Gamma$ for various schemes.
For the examples in (a), (c) and (d), the tile is $v_2$-connected.
Note that $v_m$-connectedness of a tile is different from its $v_m$
cross-sections being connected.
\par
Let us use the shorthand notation $\langle \alpha \rangle_{_0} :=
\langle \alpha \rangle_{[-\frac{1}{2},\frac{1}{2})} =
\langle \alpha + \frac{1}{2} \rangle - \frac{1}{2}$. For a
$v_m$-connected tile, we have the following simple observation:

\begin{corollary}\label{G-formula}
If the tile $\Gamma$ is $v_m$-connected, then for any ${\bf v}' \in \RR^{m-1}$
\begin{equation}
G_\Gamma({\bf v}',v_m) = \langle v_m - \lambda_\Gamma({\bf v}') \rangle_{_0}
+ \lambda_\Gamma({\bf v}'),
\end{equation}
and
$\bar{G}_\Gamma = \bar{\lambda}_\Gamma $.
\end{corollary}

\begin{proof}
If $\Gamma$ is $v_m$-connected, then $\Lambda_\Gamma({\bf v}')
= [\lambda_\Gamma({\bf v}') -\frac{1}{2},\lambda_\Gamma({\bf v}') +
\frac{1}{2})$. The first result follows from Proposition \ref{lambda-set} and
the identity $\langle \beta \rangle_{[\alpha-
\frac{1}{2},\alpha+\frac{1}{2})}= \langle \beta - \alpha \rangle_{_0} +
\alpha$ which holds for any $\alpha$ and $\beta$. The second result
is a simple consequence of the fact that the first term integrates to
zero over $v_m$.
\end{proof}

\par Before we state the following proposition, let us note that the matrix
${\bf L}^k := {\bf L}^k_m$ can be decomposed as
\begin{equation}
{\bf L}^k_m =
\left (
\begin{array}{cc}
{\bf L}^k_{m-1} & \fbox{$\begin{array}{c} {} \\ {\bf 0} \\ {} \end{array}$} \\
\fbox{$\begin{array}{ccc} & {\bf s}^\top_{m-1}[k] & \end{array}$} & 1
\end{array}
\right ),
\end{equation}
since ${\bf s}_{m}[k]$ satisfies ${\bf s}_{m}[k]
= {\bf L}_m {\bf s}_{m}[k-1] + {\bf 1}_{m}$ with
${\bf s}_{m}[0] = {\bf 0}$.

\begin{proposition}\label{prop-g-k}
Let the invariant tile $\Gamma$ be $v_m$-connected.
Define for each $k \in \ZZ$, and ${\bf v}' \in \RR^{m-1}$,
$$g_k({\bf v}') = {\bf s}_{m-1}[k] \cdot{\bf v}' + x s_m[k]
- \lambda_\Gamma({\bf L}^k_{m-1}{\bf v}' + x {\bf s}_{m-1}[k]) +
\lambda_\Gamma({\bf v}').$$
Then
\begin{equation}
\breve \rho_u[k] =
\int_{\TT^{m-1}} A_{\langle \cdot \rangle_{_0}}
(g_{k}({\bf v}'))\, \dif {\bf v}' \; + \;
\left(\breve \lambda_\Gamma,\, \UU^k \breve
\lambda_\Gamma\right)_{L^2(\TT^{m-1})}.
\end{equation}
In particular, if $m=2$ or if $\pr(\Gamma)$ is an interval of unit length,
then the second term drops.
\end{proposition}

\begin{proof}
We employ Corollary \ref{G-formula} for the evaluation of
$G_\Gamma({\bf v})$ and $G_\Gamma({\bf L}^k{\bf v}+x{\bf s}[k])$. Note
first that
$${\bf L}^k{\bf v}+x{\bf s}[k] = ({\bf L}^k_{m-1}{\bf v}' + x{\bf s}_{m-1}[k],
\, v_m + {\bf s}_{m-1}[k] \cdot{\bf v}' + x s_m[k]).$$
Therefore, we obtain
\begin{eqnarray*}
\lefteqn{
\int_\TT G_\Gamma({\bf v}) G_\Gamma({\bf L}^k{\bf v}+x{\bf s}[k]) \,\dif
v_m =} \\
& & \int_\TT \left \langle v_m {-} \lambda_\Gamma({\bf v}') \right \rangle_{0}
\left \langle v_m {+} {\bf s}_{m-1}[k] \cdot{\bf v}' {+} x s_m[k]
{-} \lambda_\Gamma({\bf L}^k_{m-1}{\bf v}' {+} x{\bf s}_{m-1}[k]) \right
\rangle_{0} \dif v_m \\
& & + \; \lambda_\Gamma({\bf v}')
\lambda_\Gamma({\bf L}^k_{m-1}{\bf v}' + x{\bf s}_{m-1}[k]),
\end{eqnarray*}
where the cross terms have dropped because
$\int_\TT \langle v_m + \varphi({\bf v'}) \rangle_{_0} \, \dif v_m = 0$
for any function $\varphi$.
The first term above
is equal to $A_{\langle \cdot \rangle_{_0}}(g_{k}({\bf v}'))$, whereas
if the second term is integrated over $\TT^{m-1}$ we find
$\left(\lambda_\Gamma,\, \UU^k \lambda_\Gamma\right)_{L^2(\TT^{m-1})}$.
The result follows since $\bar{\lambda}_\Gamma = \bar{G}_\Gamma$.
\par If $m=2$, then moreover $\lambda_\Gamma = \bar{G}_\Gamma$, so that
we have $\breve \lambda_\Gamma = 0$.
If $J:=\pr(\Gamma)$ is an interval of unit length, then it is necessarily
the case that $\Lambda_\Gamma = \RR^{m-1} \times J$.
In this case, we simply have $\bar \lambda_\Gamma = \lambda_\Gamma$ so that
$\breve \lambda_\Gamma = 0$. Hence the second term drops
in both cases.
\end{proof}

\subsection{Special case when $\pr(\Gamma) = [-\frac{1}{2},\frac{1}{2})$}
There is a class of quantization rules \cite{GrayMulti,gray2,He92}, 
for which $u_m[n] \in [-\frac{1}{2},\frac{1}{2})$ for all $n$ (for all $x$), 
so that the invariant tile
$\Gamma$ satisfies $\pr(\Gamma) = [-\frac{1}{2},\frac{1}{2})$.
These are the ``ideal'' rules that were mentioned in Section \ref{Intro},
and represent essentially the simplest possible quantization situation.
It turns out that the spectral measure $\mu$ is quite
different in its nature for $m=1$ and $m \geq 2$.

\par For $m=1$, we have $\bar{G}_\Gamma = G_\Gamma$.
Hence $\mu$ is pure-point, and Theorem \ref{pure-point} yields
$$ \mu = \sum_{n\not= 0} \frac{1}{4\pi^2 n^2} \, \delta_{nx}.$$

\par For $m\geq 2$, we simply note that $\lambda_\Gamma \equiv 0$, so that
$G_\Gamma({\bf v}) = \langle v_m \rangle_{_0}$. The fact that
$\int_\TT \langle v_m \rangle_{_0}\, \dif v_m = 0$ implies
$\bar{G}_\Gamma \equiv 0$. Hence $\mu_\mathrm{pp}=0$, i.e., $\mu$
is absolutely continuous. In addition, Proposition \ref{prop-g-k} yields
$$
\rho_u[k] =
\int_{\TT^{m-1}} A_{\langle \cdot \rangle_{_0}}
({\bf s}_{m-1}[k] \cdot{\bf v}' + x s_m[k])\, \dif {\bf v}'.
$$
For $k = 0$, the argument of the integrand is identically zero, so we obtain
$\rho_u[0] = A_{\langle \cdot \rangle_{_0}}(0) =  \frac{1}{12}$.
On the other hand, for all $k \not= 0$,  we find that $\rho_u[k] = 0$
since the integrand is of the form $A_{\langle \cdot \rangle_{_0}}(
kv_{m-1} + \alpha)$ which integrates to zero over the variable $v_{m-1}$.
Therefore,
$$ \rho_u[k] =
\left\{
\begin{array}{cl}
\frac{1}{12} & \mbox{ if } k = 0, \\
0 & \mbox{ if } k \not= 0,
\end{array}\right.
$$
and consequently $\mu$ is flat, and equal to $\frac{1}{12}$ times
Lebesgue measure on $\TT$, and the spectral density $s$ is the constant
function $s(\xi) \equiv \frac{1}{12}$.
\par
These results were previously obtained, in the case $m=1$ 
in \cite{gray2}, and in the case
$m \geq 2$ in \cite{GrayMulti,He92}.

\section{Analysis of the mean square error}

We are interested in the asymptotical behavior of $\EE(x,\phi)$
for a given \sd modulation scheme of order $m$ as the support of
$\phi$ increases and its Fourier transform $\Phi$ localizes around
zero frequency. There will be two standard choices for $\Phi$:
\begin{enumerate}
\item The ideal low-pass filter given by
$$\Phi^{\mathrm{id}}_M(\xi) :=
\chi_{[-\frac{1}{M},\frac{1}{M}]}(\xi),$$
\item The $\mathrm{sinc}^p$ family given by
$$\mathrm{Sinc}^p_M(\xi)
:= \left(\frac{\sin(\pi M \xi)}{M\sin(\pi \xi)}
e^{i \pi (M-1)\xi}\right)^p.
$$
\end{enumerate}
Note that $\mathrm{Sinc}^p_M(\xi)$ has Fourier coefficients given by
$$ \mathrm{sinc}^p_M[n] =
r_M^{(p)}[n]:=
(\underbrace{r_M*r_M*\cdots*r_M}_{p~\mathrm{times}})[n],$$
where
$r_M$ denotes the rectangular sequence
$$r_M[n]=\left\{\begin{array}{ll}
1/M, & 0\leq n<M,\\
0,& \mbox{otherwise.}\end{array}\right.$$
It is a standard fact that $ \mathrm{sinc}^p_M$ is a discrete B-spline
of degree $p-1$.

\par
We decompose the mean square error $\EE(x,\phi)$ as
$$\EE(x,\phi) = \EE_\mathrm{pp}(x,\phi)
+ \EE_\mathrm{ac}(x,\phi)$$
which correspond to the additive contributions of $\mu_\mathrm{pp}$
and $\mu_\mathrm{ac}$, respectively, in the formula
(\ref{MSE-spec}). Note that both terms are non-negative, and the
straightforward inequality
\begin{equation}
|\Phi^{\mathrm{id}}_M(\xi)| \; \lesssim_p\;
|\mathrm{Sinc}^p_{M/2}(\xi)|, \;\;\;\;\forall \xi \in \TT,
\end{equation}
implies that for any of these terms
it suffices to prove lower bounds for
the ideal low-pass filter or upper bounds for any
$\mathrm{sinc}^p$ family.

\subsection{The pure-point contribution
$\EE_\mathrm{pp}(x,\phi)$}

\par
Our first formula follows directly from plugging the expression
for $\mu_\mathrm{pp}$ given by Theorem
\ref{pure-point} in (\ref{MSE-spec}):
\begin{equation}\label{MSE-pp}
\EE_\mathrm{pp}(x,\phi) = \sum_{n \in \ZZ}
|2\sin(\pi n x)|^{2m} |\Phi(nx)|^2
\left| \widehat{\bar G_{\Gamma}}[n] \right|^2.
\end{equation}

\par
Before we carry out our analysis of this expression, let us recall
some elementary facts about Diophantine approximation. For $\alpha
\in \RR$, let $\|\alpha\|$ denote the distance of $\alpha$ to the
nearest integer, that is
$\|\alpha\| := \min(\langle \alpha \rangle,\,
\langle -\alpha \rangle)$. The number $\alpha$ is said to be
(Diophantine) of type $\eta$ if $\eta$ is the infimum of all numbers
$\sigma$ for which $$\| n \alpha \| \; \gtrsim_{\sigma,\alpha} \;
|n|^{-\sigma} \;\;\; \forall n \in \ZZ\backslash \{0\}.$$
Almost every real number (in the sense of Lebesgue measure) is
of type $1$, the smallest attainable type.

\par
The following theorem shows that for almost every $x$, if the
function $G_\Gamma$ has a sufficiently regular projection
$\bar G_\Gamma$, then the pure-point part of
the mean square error after $\mathrm{sinc}^{m+1}_M$ filtering
decays at least as fast as $M^{-2m-2}$.

\begin{theorem}\label{MSEppDioph}
Let $x$ be Diophantine of type $\eta$, and $\alpha$ and $\beta$ be two
real numbers satisfying $0 \leq \alpha \leq 1$ and
$\beta > (1-\frac{\alpha}{2})\eta$. If the invariant tile
$\Gamma = \Gamma_x$ of an $m$'th order \sd modulator with input $x$
satisfies
$$\left| \widehat{\bar G_\Gamma}[n] \right| \;
\lesssim \; |n|^{-\beta}$$
for all $n \in \ZZ\backslash\{0\}$, then
\begin{equation}
\EE_\mathrm{pp}(x,\mathrm{sinc}^{m+1}_M) \;
\lesssim_{x,m,\alpha,\beta} \; M^{-2m-2+\alpha}
\end{equation}
for all $M$.
\end{theorem}

\begin{proof}
Formula (\ref{MSE-pp}) reads
\begin{equation}\label{MSE-pp-sinc}
\EE_\mathrm{pp}(x,\mathrm{sinc}^{m+1}_M)
= \frac{2^{2m}}{M^{2m+2}} \sum_{n \in \ZZ\backslash\{0\}}
\frac{\sin^{2m+2}(\pi M nx)}{\sin^2(\pi n x)}
\left| \widehat{\bar G_\Gamma}[n] \right|^2.
\end{equation}
\par
Given the decay of $| \widehat{\bar G_\Gamma}[n]|$
and the simple fact $|\sin(\pi \theta)| \asymp \|\theta\|$,
it suffices to show that
\begin{equation*}
\sum_{n=1}^\infty \frac{\|Mnx\|^{2m+2}}{\|nx\|^2}\frac{1}{n^{2\beta}}
\; \lesssim_{x,\beta} \;M^{\alpha}.
\end{equation*}
Note that $\|Mnx\| \leq \min(1,M\|nx\|)$. Hence we have
\begin{eqnarray*}
\sum_{n=1}^\infty \frac{\|Mnx\|^{2m+2}}{\|nx\|^2}\frac{1}{n^{2\beta}}
& \leq  &
\sum_{n=1}^\infty \frac{\|Mnx\|^\alpha}{\|nx\|^\alpha}
\frac{1}{\|nx\|^{2-\alpha}n^{2\beta}} \\
& \leq &
M^\alpha \sum_{n=1}^\infty \frac{1}{n^{2\beta}\|nx\|^{2-\alpha}}.
\end{eqnarray*}
Since $2-\alpha \geq 1$, it suffices to show the convergence of the sum
$$
\sum_{n=1}^\infty \frac{1}{n^{2\beta/(2-\alpha)} \|nx\|}.
$$
\par
Let $\lambda>0$ be defined by $\beta = (1-\frac{\alpha}{2})(\eta + \lambda)$.
Now, summation by parts shows that
\begin{equation}\label{Abel-nx}
\sum_{n=1}^\infty \frac{1}{n^{\eta+\lambda} \|nx\|}
\;\lesssim_{\eta,\lambda} \; \sum_{n=1}^\infty \frac{1}{n^{\eta+\lambda+1}}
\left( \sum_{k=1}^n \frac{1}{\|kx\|}\right )
\end{equation}
and furthermore it is well-known \cite[Ex. 3.11]{kuipers} that
$$
 \sum_{k=1}^n \frac{1}{\|kx\|} \; \lesssim_{x,\sigma} \; n^{\sigma}
$$
for any $\sigma > \eta$. Choosing $\eta + \lambda > \sigma > \eta$,
we obtain the convergence of (\ref{Abel-nx}) with a sum depending
on $x$ and $\lambda$. Since $\lambda$ depends on $x$, $\alpha$ and $\beta$,
the result of the theorem follows.
\end{proof}

\par
{\bf Note:}
The Diophantine condition on $x$ can be removed if $\bar G_\Gamma$
is a trigonometric polynomial. In this case, (\ref{MSE-pp-sinc})
reduces to a finite sum, and therefore it is always convergent.

\par
On the other hand,
our next result shows that if $\bar G_\Gamma$ does not have
enough regularity in a certain sense as specified in the following theorem,
then this is the best one can get in the sense that there is an everywhere
dense set of exceptional values of $x$ for which the exponent of the
error decay rate is never better than $2m$, even for the ideal low pass
filter.

\begin{theorem}
Given a \sd modulator of order $m$,
let $\phi_M$, $M=1,2,\dots$, be a sequence of
averaging filters such that
$|\Phi_M(\xi)| \geq c_1$ on the interval $|\xi| \leq c_2/M$, where
$c_1$ and $c_2$ are positive constants that do not depend on $M$.
There exists a dense set $E$ 
of irrational numbers with the following property:
For any $x \in E$, if there exist positive constants $\beta_x$
and $C_x$ such that the invariant tile
$\Gamma = \Gamma_x$ satisfies
$$\left| \widehat{\bar G_\Gamma}[n] \right| \geq C_x |n|^{-\beta_x}$$
for all but finitely many $n \in \ZZ$, then for all $\delta > 0$,
\begin{equation}\label{MSE-lower}
\limsup_{M \to \infty} \;
\EE_\mathrm{pp}(x,\phi_M) \, M^{2m+\delta} = \infty.
\end{equation}
\end{theorem}

\begin{proof}
It suffices to find, for any
open interval $J$, a point $x \in J$ with the property
(\ref{MSE-lower}) for all $\delta > 0$. Given an open interval $J$,
let $x_0 \in J$ be a dyadic rational. Let $l = \max(b,d)$ for the
minimum $b$ and $d$
such that $b!-1$ is an upper bound for the length of
the binary expansion of $x_0$ and $2^{-d!+1}$ is a lower
bound for the distance of $x_0$ to the boundary of $J$. Set
$$x = x_0 + \sum_{k \geq l} 2^{-k!}.$$
Then clearly $x \in J$. It is also a standard fact that $x$ is
irrational, in fact $x$ is a Liouville number.
\par
Note that for $q \geq l$, we have
$$ \langle 2^{q!}x \rangle = \sum_{k=q+1}^\infty 2^{-k! + q!} =
2^{-q\cdot q!} + \sum_{k=q+2}^\infty 2^{-k! + q!}.
$$
For $q=1,2,\dots$, let $n_q = 2^{q!}$
and  $M_q = \lfloor 2^{q\cdot q! - 1} c_2 \rfloor$. We have
\begin{equation}\label{nq-x}
\frac{c_2}{4M_q} < 2^{-q\cdot q!} < \| n_q x \| < 2^{-q\cdot q! +1 }
< \frac{c_2}{M_q}.
\end{equation}
The right side of this chain of inequalities implies
$|\Phi_{M_q}(n_q x)| \geq c_1$ by our assumption on $\{\phi_M\}$.
On the other hand, the left side implies
$|2 \sin(\pi n_q x)| \geq 4 \|n_q x \| > c_2/M_q$. Therefore
\begin{eqnarray}
 \EE_\mathrm{pp} (x,\phi_{M_q})
& > &
|2\sin(\pi n_q x)|^{2m} |\Phi(n_qx)|^2
\left| \widehat{\bar G_\Gamma}[n_q] \right|^2 \nonumber \\
& > &
C_{x,m} M_q^{-2m} M_q^{-\beta_x/q}
\end{eqnarray}
where $C_{x,m} = C_x c_1^2 c_2^{2m}$.
The result of the theorem follows by letting $q \to \infty$ and therefore
exhibiting the subsequence $M_q$ for which (\ref{MSE-lower}) holds for
any $\delta > 0$.
\end{proof}

\subsection{The absolutely continuous contribution
$\EE_\mathrm{ac}(x,\phi)$}

Let us denote by $s$ the Radon-Nikodym derivative of the absolutely
continuous spectral measure $\mu_\mathrm{ac}$, i.e.,
$\dif \mu_\mathrm{ac} = s(\xi) \dif \xi$. A priori, we know that
$s \in L^1(\TT)$, which is somewhat weak for what we would like to
achieve. Our first theorem concerns the decay rate of
$$
\EE_\mathrm{ac}(x,\mathrm{sinc}^{m+1}_M)
= \int_\TT |2\sin(\pi \xi)|^{2m} |\mathrm{Sinc}^{m+1}_M(\xi)|^{2}
s(\xi) \, \dif \xi
$$
if it is known that $s$ belongs to an $L^p$ space.

\begin{theorem}\label{MSEacLpcase}
If the measure $\mu_\mathrm{ac}$
has density $s \in L^p(\TT)$ for some $1 \leq p \leq \infty$, then
\begin{equation}
\EE_\mathrm{ac}(x,\mathrm{sinc}^{m+1}_M) \;
\lesssim_{m,p} \;
 \| s \|_{L^p(\TT)} \, M^{-(2m+1-1/p)}.
\end{equation}
\end{theorem}

\begin{proof}
Let ${p'}$ be the dual index of $p$, i.e., $1/p + 1/{p'} = 1$.
Note that
\begin{eqnarray}
|2\sin(\pi \xi)|^{2m}  |\mathrm{Sinc}^{m+1}_M(\xi)|^{2}
& =_{\;\;\,} &
\label{sinc-to-sinc}
|2\sin(\pi M \xi)|^{2m} |\mathrm{Sinc}_M(\xi)|^{2} M^{-2m}  \\
& \lesssim_m &  |\mathrm{Sinc}^2_M(\xi)|\, M^{-2m} ,
\end{eqnarray}
so that H\"older's inequality yields
$$
\EE_\mathrm{ac}(x,\mathrm{sinc}^{m+1}_M) \;
\lesssim_m  \big \| s \big \|_{L^p(\TT)}
\left \| \mathrm{Sinc}^2_M \right \|_{L^{p'}(\TT)}  M^{-2m} .
$$
Furthermore, the simple bound $|\mathrm{Sinc}_M(\xi)| \leq
\min\!\left(1,(2M|\xi|)^{-1}\right)$ implies
\begin{equation}
\left \| \mathrm{Sinc}^2_M \right \|_{L^{p'}(\TT)}
\; \lesssim_{p'} \; M^{-1/{p'}},
\end{equation}
hence the theorem follows.
\end{proof}

\par
On the other hand, it turns out that if $s$ is continuous at $0$,
then one can calculate the exact asymptotics of
$\EE_\mathrm{ac}(x,\mathrm{sinc}^{m+1}_M)$ without additional
assumptions.

\begin{theorem}\label{MSEacCtscase}
If the spectral density $s$ is continuous at $0$, then
\begin{equation}
\EE_\mathrm{ac}(x,\mathrm{sinc}^{m+1}_M) \;
= \; \binom{2m}{m}\,s(0)\, M^{-2m-1} + o(M^{-2m-1}).
\end{equation}
\end{theorem}
\begin{proof}
The proof has two parts. First part is the easy calculation
\begin{equation}\label{2m-choose-m}
\int_\TT |2\sin(\pi \xi)|^{2m} |\mathrm{Sinc}^{m+1}_M(\xi)|^{2}
\, \dif \xi
= \; \binom{2m}{m}\, M^{-2m-1}.
\end{equation}
To see this, note that (\ref{sinc-to-sinc}) and the definition of
$\mathrm{Sinc}_M(\xi)$ imply
$$
|2\sin(\pi \xi)|^{2m} |\mathrm{Sinc}^{m+1}_M(\xi)|^{2}
= \left(\frac{e^{i \pi M \xi}- e^{-i \pi M \xi}}{i}\right)^{2m}
\sum_{k=0}^{M-1} \sum_{j=0}^{M-1} e^{2\pi i (k-j) \xi} M^{-2m-2}.
$$
The right hand side is the product of two trigonometric polynomials;
the first polynomial has frequencies only at integer multiples of $2\pi M$ and
the second polynomial has frequencies between $-2\pi (M-1)$ and $2\pi (M-1)$.
The zero frequency term of the product is therefore given only by the
product of the corresponding zero frequency terms, which is equal to
$$ \binom{2m}{m} (-1)^m i^{-2m} \left(\sum_{k=0}^{M-1} 1\right) M^{-2m-2}
=  \binom{2m}{m}\, M^{-2m-1},$$
hence the result.
\par
The second part of the proof concerns the residual term
$$
\left|
\int_\TT \big |2\sin(\pi \xi) \big |^{2m}
\big |\mathrm{Sinc}^{m+1}_M(\xi)\big|^2
\big (s(\xi)-s(0) \big ) \, \dif \xi \right |,
$$
which is bounded, using (\ref{sinc-to-sinc}), by
$$
2^{2m} M^{-2m} \int_\TT \mathrm{Sinc}^2_M(\xi)
|s(\xi)-s(0)| \, \dif \xi
= 2^{2m} M^{-2m-1} \int_\TT K_{M-1}(\xi) |s(\xi)-s(0)| \, \dif \xi,
$$
where
$$ K_{M-1}(\xi) = \frac{1}{M} \left(\frac{\sin(\pi M \xi)}{\sin(\pi \xi)}
\right)^2$$
is the Fej\'{e}r kernel. The limit
$$ \lim_{M\to \infty} \; \int_\TT K_{M-1}(\xi) |s(\xi)-s(0)| \, \dif \xi $$
is the Ces\`{a}ro sum of the Fourier series of the function
$f(t) = |s(-t)-s(0)|$ evaluated at $t = 0$. Since
$f$ is continuous at $0$, the Ces\`{a}ro sum converges to $f(0) = 0$, and
therefore the limit is $0$. This concludes the proof.
\end{proof}
{\bf Notes:}
\begin{enumerate}
\item
A similar calculation shows that for the ideal
filter $\Phi^\mathrm{id}_M$, the error has the asymptotics given by
\begin{equation}
\EE_\mathrm{ac}(x,\Phi^\mathrm{id}_M) \;
= \; \frac{(2\pi)^{2m+1}}{m+1/2}\,s(0)\, M^{-2m-1} + O(M^{-2m-3})
\end{equation}
again assuming that $s$ is continuous at $0$.
\item The value of $s(0)$ is equal to the sum of its
Fourier coefficients $\breve \rho_u[k]$.
\end{enumerate}

\section{Estimates for
second order schemes with $v_2$-connected invariant tiles}

\par
Second order \sd modulators with $v_2$-connected invariant tiles
are interesting because the value of $x$ and the function
$\lambda_\Gamma = \bar G_\Gamma$ completely describe the MSE
behavior via the theorems we have stated in the previous sections.
In particular, (\ref{prop-g-k}) provides us with the formula
\begin{eqnarray}
\breve \rho_u[k]
& = &
\int_\TT \textstyle
A_{\langle \cdot \rangle_{_0}}\!\!
\left ( kv_1 + \frac{k(k+1)}{2}x - \lambda_\Gamma(v_1+kx) + \lambda_\Gamma(v_1)
\right )\, \dif v_1 \nonumber \\
& = &
\int_\TT \textstyle
A_{\langle \cdot \rangle_{_0}}\!\!
\left ( kv - \lambda_\Gamma(v-\frac{x}{2}+k\frac{x}{2})
+ \lambda_\Gamma(v-\frac{x}{2}-k\frac{x}{2})
\right )\, \dif v,
\end{eqnarray}
where we have used the change of variable $v = v_1 + (k+1)x/2$ to
obtain the second representation.
\par
By Riemann-Lebesgue lemma, we already know that $\breve \rho_u[k]$ must
converge to zero as $|k| \to \infty$ since $\breve \rho_u[k] = \hat s[k]$,
where $s \in L^1(\TT)$ is the spectral density. However, we
would like to quantify the rate of decay in $|k|$ as
this would then allow us to draw conclusions about $s$.
Intuitively speaking,
it is not hard to see from this formula that the smoother $\lambda_\Gamma$ is,
the faster $\breve \rho_u[k]$ must decay in $|k|$ as $|k| \to \infty$,
since $A_{\langle \cdot \rangle_{_0}}$ is a zero mean function on $\TT$.
Our objective in this section is to study this relation rigorously.
\par
Let $\mathrm{BV}(\TT)$ denote the space of functions on $\TT$
that have bounded variation, where $\| \cdot \|_{TV}$ denotes the
total variation semi-norm, and let $\mathrm{A}(\TT)$ denote the space of
functions on $\TT$ with absolutely convergent Fourier series with
the norm $\| f \|_{\mathrm{A}(\TT)}$ given by $\sum |\hat f[n]|$.
We have the following lemma, whose proof is given in the Appendix.

\begin{lemma}\label{lemma-decay}
Let $f \in {\mathrm{A}(\TT)}$
and $\varphi$ be two real valued functions on $\TT$,  where
$f$ has zero mean. Consider the integrals
\begin{equation}\label{c[k]}
c[k] = \int_\TT f(kv + \varphi(v)) \,\dif v.
\end{equation}
The following bounds hold:
\begin{enumerate}
\item If $\varphi \in \mathrm{BV}(\TT)$, then for all $k \in \ZZ\backslash
\{0\}$,
\begin{equation}\label{c[k]bound1}
\big|c[k]\big|
\leq \frac{1}{|k|}  \| f \|_{\mathrm{A}(\TT)} \| \varphi \|_{TV}.
\end{equation}
\item
If $\varphi$ is differentiable almost everywhere and
$\varphi' \in \mathrm{BV}(\TT)$, then for all $k \in \ZZ\backslash
\{0\}$,
\begin{equation}\label{c[k]bound2}
\big |c[k] \big | \leq
\frac{1}{k^2} \left(\frac{1}{\sqrt{12}} \,\|f\|_{L^2(\TT)}\|\varphi'\|_{TV}
+ \|f\|_{L^\infty(\TT)} \|\varphi'\|^2_{L^2(\TT)} \right).
\end{equation}
\end{enumerate}
\end{lemma}

\begin{corollary}
Let $x$ be given and $\Gamma$ be the invariant tile corresponding to
a second order \sd modulator. Then we have the following:
\begin{enumerate}
\item
If the midpoint function $\lambda_\Gamma$ has
bounded variation on $\TT$, then
\begin{equation} \label{rho1}
\big |\breve \rho_u[k] \big|
\leq \frac{1}{6|k|} \big \| \lambda_\Gamma \big \|_{TV}.
\end{equation}
Consequently, one has
\begin{equation}\label{MSEac1}
\EE_\mathrm{ac}(x,\mathrm{sinc}^{3}_M) \;
\lesssim_{x,\epsilon} \; M^{-5+\epsilon}
\end{equation}
for any $\epsilon > 0$. If the type $\eta$ of $x$ is strictly less than
$2$, then
\begin{equation}\label{MSEpp1}
\EE_\mathrm{pp}(x,\mathrm{sinc}^{3}_M) \;
\lesssim_{x,\delta} \; M^{-5-\delta}
\end{equation}
for any $0 \leq \delta < (2-\eta)/\eta$.
\item If the midpoint function $\lambda_\Gamma$ has a derivative
that has bounded variation on $\TT$, then
\begin{equation}\label{rho2}
\big | \breve \rho_u[k] \big | \leq
\frac{1}{k^2} \left(\frac{1}{12\sqrt{15}}
\,\|\lambda'_\Gamma\|_{TV}
+ \frac{1}{3} \|\lambda'_\Gamma\|^2_{L^2(\TT)} \right).
\end{equation}
In particular, the spectral density
$s$ is continuous. Consequently, one has
\begin{equation}
\EE_\mathrm{ac}(x,\mathrm{sinc}^{3}_M) \;
= \; 6\,s(0) M^{-5} + o(M^{-5}),
\end{equation}
where
\begin{equation}\label{s0-bound}
s(0) \leq \frac{1}{12} + \frac{\pi^2}{3}\left(\frac{1}{12\sqrt{15}}
\| \lambda'_\Gamma \|_{TV}
+ \frac{1}{3} \| \lambda'_\Gamma \|^2_{L^2(\TT)} \right).
\end{equation}
If the type $\eta$ of $x$ is strictly less than $4$, then
\begin{equation}
\EE_\mathrm{pp}(x,\mathrm{sinc}^{3}_M) \;
\lesssim_{x,\delta} \; M^{-5-\delta}
\end{equation}
for any $0 \leq \delta < \min\!\big(1,(4-\eta)/\eta\big)$.
\end{enumerate}
\end{corollary}

\begin{proof}
\par
Let
$$f(v) := A_{\langle \cdot \rangle_{_0}}(v)
= \sum_{n \not= 0} \frac{1}{4\pi^2 n^2} \, e^{2\pi inv}.$$
For each $k$, define
$$ \textstyle \varphi_k(v) :=  - \lambda_\Gamma(v-\frac{x}{2}+k\frac{x}{2})
+ \lambda_\Gamma(v-\frac{x}{2}-k\frac{x}{2}).$$
For these functions, we have the following exact formulas and bounds:
\begin{eqnarray}
\|f \|_{\mathrm{A}(\TT)} & = & \frac{1}{12} \label{eqfor1}\\
\|f \|_{L^\infty(\TT)} & = & \frac{1}{12}  \label{eqfor2}\\
\|f \|_{L^2(\TT)} & = &   \frac{1}{12\sqrt{5}} \label{eqfor3}\\
\|\varphi_k\|_{TV} & \leq & 2\, \|\lambda_\Gamma\|_{TV} \label{ineqfor1} \\
\|\varphi'_k\|_{TV} & \leq & 2\, \|\lambda'_\Gamma\|_{TV}  \label{ineqfor2}\\
\|\varphi'_k\|^2_{L^2(\TT)} & \leq & 4\, \|\lambda'_\Gamma\|^2_{L^2(\TT)}.
 \label{ineqfor3}
\end{eqnarray}

\begin{enumerate}
\item
In this case we only know that $\lambda_\Gamma$ is of bounded variation.
\par
The decay estimate
(\ref{rho1}) simply follows from the bound (\ref{c[k]bound1}) coupled
with (\ref{eqfor1}) and (\ref{ineqfor1}).
\par
Given that the Fourier coefficients $\breve \rho_u[k]=\hat s[k]$
decay like $1/k$,
it follows from Riesz-Thorin interpolation theorem that
the spectral density $s\in L^p(\TT)$
for any $p < \infty$. Therefore Theorem \ref{MSEacLpcase}
implies, with $m=2$ and $\epsilon = 1/p$, the bound (\ref{MSEac1}).
\par
For the pure-point estimate, we use Theorem  \ref{MSEppDioph}
with $\beta = 1$ and $m=2$. If we define $\delta = 1-\alpha$, where
$\alpha$ is as defined in Theorem  \ref{MSEppDioph}, then the result
follows as stated.

\item
In this case we know that $\lambda_\Gamma$ has a derivative that
is of bounded variation.
\par
The decay estimate (\ref{rho2}) follows from the
bound (\ref{c[k]bound2}) coupled with
(\ref{eqfor2}), (\ref{eqfor3}), (\ref{ineqfor2}) and (\ref{ineqfor3}).
\par
Since $\breve \rho_u$ is summable, it follows that $s$ is
continuous. We therefore apply Theorem \ref{MSEacCtscase} to compute the
exact asymptotics of $\EE_\mathrm{ac}(x,\mathrm{sinc}^{3}_M)$.
In this case, the nonnegative number $s(0)$ will be bounded by
$\sum |\breve \rho_u[k]|$. We simply, add up the bounds
given by (\ref{rho2}), including the trivial case
$|\breve \rho_u[0]| \leq \|f\|_{L^\infty(\TT)}$. This computation
yields the bound (\ref{s0-bound}).
\par
For the pure-point estimate, we again use Theorem  \ref{MSEppDioph}, but
now with $\beta = 2$. We define $\delta = 1-\alpha$, where
$\alpha$ is as defined in Theorem  \ref{MSEppDioph}, and note that the
condition $\alpha \leq 1$ must be imposed, which was automatically
satisfied in part (1). Then the result follows as stated.
\end{enumerate}
\end{proof}

\section{Further remarks}

In this paper, we have covered only a portion of the mathematical
problems that concern \sd quantization. We believe that 
the following currently 
unresolved problems are interesting both from 
the dynamical systems standpoint and the engineering perspective:

\par 1. \textit{Which maps $\MM$ are stable?} Satisfactory answers of 
this question would include non-trivial sufficient conditions in terms
of the quantization rule $Q$, or in terms of the partition $\Pi_x$ and
the quantization levels $\{d_i\}$.

\par 2. \textit{Which stable maps $\MM$ yield single invariant tiles?}
One can include in this the case when $\Gamma$ is composed of
tiles each of which is invariant under $\MM$. In principle, 
each of these invariant tiles would represent a different ``mode
of operation''.

\par 3. \textit{What is an appropriate generalization of
our spectral analysis of mean square error when $\Gamma$ is composed of 
more than one tile?}

\par 4. \textit{Given the quantization rule,
what can be said about the geometric regularity of $\Gamma$?}
We used two types of geometric information about $\Gamma$ in deriving our
analytical results on the mean square error asymptotics. The
first type concerned ``shape'' (such as $v_m$-connectedness),
and the second concerned ``regularity'' (such as the decay of Fourier
coefficients of $\bar{G}_\Gamma$). At this stage, the relation between the
quantization rule and these two issues is highly unclear, although
we have partial understanding in some cases. Even for ``linear'' rules,
there seems to be a wide range of possibilities.

\par 5. \textit{What are the universal principles behind tiling?}
Tiling invariant sets are found even when $x$ is rational. 
In addition, trajectories seem to remain within exact tiles, and not
just tiles ``up to sets of measure zero''.

\renewcommand{\thesubsection}{\Alph{section}.\arabic{subsection}}
\renewcommand{\theequation}{\Alph{section}.\arabic{equation}}
\renewcommand{\thetheorem}{\Alph{subsection}.\arabic{theorem}}

\setcounter{section}{1}
\section*{Appendix A. On the spectral theory of the map $\mathscr{L}$}
\setcounter{equation}{0}

In this section, we will review some basic facts about the spectral
theory of the map $\mathscr{L} = \mathscr{L}_x$ on $\TT^m$, where
$\mathscr{L}_x {\bf v}={\bf L}{\bf v}+x{\bf 1}$,
and $x$ is an irrational number. Most of what follows below can
be derived or generalized from Anzai's work on ergodic skew product
transformation \cite{Anzai}.

\subsection*{The eigenfunctions of $\UU_\LL$}
We start by showing that the set of all eigenfunctions
of $\UU=\UU_\mathscr{L}$ is precisely given by the collection
of complex exponentials $f_n$, where
$$f_n({\bf v}) = e^{2\pi in v_1}, \;\;\;\;n \in \ZZ.$$
To see this, let $f\in L^2(\TT^m)$ be an eigenfunction of $\;\UU$
with eigenvalue $\lambda$.
Since $\UU$ is unitary, $|\lambda|=1$. Consider
the Fourier series expansion of $f$ given by
$$ f({\bf v}) = \sum_{{\bf n}\in\ZZ^m}
c[{\bf n}] e^{2\pi i {\bf n}\cdot{\bf v}}.$$
Since $f = \frac{1}{\lambda} \big(\UU f\big)$,
we have the relation
\begin{eqnarray*}
\sum_{{\bf n}\in\ZZ^m} c[{\bf n}] e^{2 \pi i {\bf n}\cdot{\bf v}}
& = & \frac{1}{\lambda}
\sum_{{\bf n}\in\ZZ^m} c[{\bf n}] e^{2\pi i x {\bf n}\cdot{\bf 1}}
e^{2\pi i {\bf n}\cdot({\bf Lv})} \\
& = & \frac{1}{\lambda}
\sum_{{\bf n}\in\ZZ^m} c[{\bf K}{\bf n}]
e^{2\pi i x({\bf K}{\bf n})\cdot{\bf 1}}
e^{2\pi i {\bf n}\cdot{\bf v}},
\end{eqnarray*}
where  ${\bf K}=({\bf L}^{-1})^\top$. Comparing the coefficients,
we obtain the equality
$$\big |c[{\bf n}]\big |= \big |c[{\bf K n}]\big |,\;\;\;\;
\forall{\bf n}\in\ZZ^m.$$
Since $f \in L^2(\TT^m)$, we can conclude that $c[{\bf n}]=0$
for any ${\bf n}$ that is not
preserved under ${\bf K}^j$ for some positive integer $j$, for otherwise
we would have the infinite sequence of coefficients
$c[{\bf n}], c[{\bf Kn}], c[{\bf K}^2{\bf n}], \dots$
of equal and strictly positive magnitude. In fact, this conclusion is valid
even for $f \in L^1(\TT^m)$ since in this case it is also true that
$|{\bf K}^j{\bf n}| \to \infty$ as $j \to \infty$, resulting in a
violation of the Riemann-Lebesgue lemma.
\par
On the other hand, it is a simple exercise to show
that the only vectors that satisfy ${\bf n} = {\bf K}^{j}{\bf n}$ for
some power $j \geq 1$ are those
of the form ${\bf n} = (n_1,0,\dots,0)$. Hence, any eigenfunction
of $\;\UU$ depends only on the first variable $v_1$.
On the first coordinate $v_1$ of ${\bf v}$,
$\mathscr{L}$ reduces to the irrational rotation by $x$, and hence as
it is well-known, these eigenfunctions are nothing but the given
complex exponentials $\{f_n\}_{n \in \ZZ}$. According to the spectral
theorem, these eigenfunctions span the subspace $\mathscr{H}_\mathrm{pp}$
of $L^2(\TT^m)$.

\subsection*{The absolutely continuous spectrum}
We shall next show that continuous part of the spectrum is in fact
absolutely continuous. This is in fact a consequence of the fact that
there exists an orthonormal basis $\{\psi_{j,k} :
j \in \ZZ, k \in \NN \}$ of $\HH_\mathrm{pp}^\perp$ with the
property\footnote{I.e., $\LL$ has
{\em countable Lebesgue spectrum} on $\HH_\mathrm{pp}^\perp$.} that
$\UU \psi_{j,k} = \psi_{j+1,k}$ for all $j$ and $k$. First we will
construct such a basis, and then we shall prove the statement on the
absolute continuity.
\par
From the discussion above, we know that the complex exponentials
$$f_{\bf n}({\bf v}) = e^{2\pi i {\bf n}\cdot{\bf v} },\;\;\;\;
{\bf n} \in \ZZ^m\backslash\left(\ZZ\times\{0\}^{m-1}\right),$$
form an orthonormal complete set in  $\HH_\mathrm{pp}^\perp$.
Note also that
$$\UU_\LL f_{\bf n} = e^{2\pi i {\bf n}\cdot {\bf 1}}
f_{{\bf L}^\top{\bf n}}.$$
Therefore we consider the orbit of each ${\bf n} \in \ZZ^m$
under ${\bf L}^\top$, given by
$$\mathscr{O}({\bf n})
= \left \{ \left({\bf L}^\top\right)^j{\bf n}\right\}_{j \in \ZZ }.$$
It is easy to see that each ${\bf n}
\in \ZZ\times\{0\}^{m-1}$ is a fixed point of ${\bf L}^\top$ and
every other ${\bf n}$ is such that the orbit is an infinite sequence of
distinct points in $\ZZ^m\backslash\left(\ZZ\times\{0\}^{m-1}\right).$
Divide $\ZZ^m\backslash\left(\ZZ\times\{0\}^{m-1}\right)$
into equivalence classes of orbits $\mathscr{O}({\bf n}_k)$,
$k \in \NN$ (note that distinct orbits do not
intersect because ${\bf L}^\top$ is invertible), and define
$$ \psi_{0,k} = f_{{\bf n}_k},\;\;\;\;\;\;\;\;\;
\psi_{j,k} = \UU_\LL^j \psi_{0,k}, \;\;\;\;\;j \in \ZZ,\;k\in \NN.$$
Each $\psi_{j,k}$ is equal to some complex exponential $f_{\bf n}$
multiplied by a complex number of unit magnitude.
The collection of $\psi_{j,k}$ is distinct, and all frequencies
${\bf n} \in \ZZ^m\backslash\left(\ZZ\times\{0\}^{m-1}\right)$ appear, hence
$\{\psi_{j,k}\}_{j \in \ZZ, k\in \NN}$ form an orthonormal basis of
$\HH_\mathrm{pp}^\perp$ with the property that
$\UU \psi_{j,k} = \psi_{j+1,k}$.

\par
Let us show that every spectral measure is absolutely continuous on
$\HH_\mathrm{pp}^\perp$. Let $g$ and $h$ be arbitrary functions with
representations $g = \sum a_{j,k} \psi_{j,k}$ and
$h = \sum b_{j,k} \psi_{j,k}$. Let the functions $A_k$ and $B_k$ be
defined on $\TT$ for each $k$ with Fourier coefficients
$(a_{j,k})_{j \in \ZZ}$ and $(b_{j,k})_{j \in \ZZ}$, respectively.
From orthogonality, we have
$$\|g\|^2 = \sum_k \int_\TT |A_k(\xi)|^2 \dif \xi < \infty,$$
and similarly for $h$ and $B_k$.
\par
Now, we have
$\UU^n h = \sum b_{j,k} \psi_{j+n,k}$, so that
\begin{eqnarray*}
(g, \UU^n h) & = & \sum_k \sum_j a_{j+n,k} \bar{b}_{j,k} \\
& = & \sum_k \int_\TT \sum_j a_{j+n,k} e^{-2\pi i j \xi} \bar{B}_k(\xi)\dif
\xi \\
& = & \int_\TT e^{2\pi i n \xi} \left ( \sum_k A_k(\xi)\bar{B}_k(\xi)
\right ) \dif \xi
\end{eqnarray*}
Here, $L^1$ norm of the function $\sum_k A_k(\xi)\bar{B}_k(\xi) $ is
bounded by $\|g\|_{L^2}\|h\|_{L^2}$ because of Cauchy-Schwarz inequality,
hence finite.
Therefore, we have that the measure $\nu_{g,h}$ defined by the inner
products $(g, \UU^n h)$ is absolutely continuous.

\setcounter{section}{2}
\section*{Appendix B. Proof of Lemma \ref{lemma-decay}}
\setcounter{equation}{0}

Let us start by writing $f(t) = \sum\limits_{n \not= 0} \hat f[n]
e^{2\pi i n t}$ so that we have
\begin{equation}\label{equiv-c[k]}
c[k] = \sum_{n \not= 0} \hat f[n] \int_\TT
e^{2\pi i n(kv+\varphi(v))} \dif v,
\end{equation}
where we have changed the order of summation
and integration. Applying integration by parts we obtain
\begin{eqnarray}
\int_\TT e^{2\pi i n\varphi(v)} e^{2\pi i nkv} \dif v
& = &
-\frac{1}{2\pi ink} \int_\TT
 e^{2\pi i nkv} \dif \!\left [e^{2\pi in\varphi(v)} \right] \nonumber \\
& = & \label{ineq1}
-\frac{1}{k} \int_\TT e^{2\pi i nkv} e^{2\pi i n\varphi(v)}  \dif \varphi(v).
\end{eqnarray}

\par Part (1).
For the integral in (\ref{ineq1}), we use the
bound
$$ \left |
-\frac{1}{k} \int_\TT e^{2\pi i nkv} e^{2\pi i n\varphi(v)}  \dif \varphi(v)
\right | \leq \frac{1}{|k|}  \int_\TT | \dif \varphi(v) | =
\frac{1}{|k|} \,\| \varphi \|_{TV},$$
and we simply get
$$ \big |c[k] \big |
\leq \sum_{n \not= 0}  \frac{1}{|k|}\, \| \varphi \|_{TV} |\hat f[n]|
\leq \frac{1}{|k|}\, \| \varphi \|_{TV} \| f \|_{\mathrm{A}(\TT)}.$$

\par Part (2). Let $\varphi$ be differentiable and
$\varphi' \in \mathrm{BV}(\TT)$. Substitute $\dif \varphi(v) =
\varphi'(v) \dif v$ and apply
another integration by parts to (\ref{ineq1}) obtain
$$
-\frac{1}{k} \int_\TT  \varphi'(v) e^{2\pi i n\varphi(v)}
e^{2\pi i nkv} \dif v
= \frac{1}{k (2\pi i nk)} \int_{\TT} e^{2\pi i n k v}
\dif \!\left [\varphi'(v) e^{2\pi i n \varphi(v)} \right ].
$$
Now,
$$
\dif \! \left [\varphi'(v) e^{2\pi i n \varphi(v)} \right ]
= e^{2\pi i n \varphi(v)} \dif \varphi'(v) +
(\varphi'(v))^2 (2\pi i n) e^{2\pi i n \varphi(v)} \,\dif v,$$
so that substituting the above two formulas together with (\ref{ineq1})
in (\ref{equiv-c[k]}), we get
\begin{equation}\label{c[k]-formula}
c[k] = \frac{1}{k^2} \left( \sum_{n \not= 0} \frac{\hat f[n]}{2\pi in}
\int_{\TT} e^{2\pi in (kv + \varphi(v))} \,\dif \varphi'(v) +
\sum_{n\not= 0} \hat f[n]
\int_{\TT} (\varphi'(v))^2 e^{2\pi n(kv+\varphi(v))} \,\dif v
\right)
\end{equation}
For the first part of this sum we use
$$\left |\int_{\TT} e^{2\pi in (kv + \varphi(v))} \,\dif \varphi'(v)
\right | \leq \int_{\TT} |\dif \varphi'(v)| =  \|\varphi'\|_{TV},$$
and
$$\sum_{n \not= 0} \frac{|\hat f[n]|}{2\pi |n|}
\leq
\left(\sum_{n\neq0}\frac{1}{(2\pi n)^2} \right)^{1/2}
\left(\sum_{n}|\hat f[n]|^2 \right)^{1/2} =
\frac{1}{\sqrt{12}}\, \|f\|_{L^2(\TT)},
$$
so that
$$ \left| \sum_{n \not= 0} \frac{\hat f[n]}{2\pi in}
\int_{\TT} e^{2\pi in (kv + \varphi(v))} \,\dif \varphi'(v) \right|
\leq   \frac{1}{\sqrt{12}} \,\|f\|_{L^2(\TT)}\|\varphi'\|_{TV}.
$$
On the other hand, the second term reduces to
\begin{eqnarray*}
\sum_{n \not= 0} \hat f[n]
\int_{\TT} (\varphi'(v))^2 e^{2\pi n(kv+\varphi(v))} \,\dif v
& = & \int_{\TT} (\varphi'(v))^2
\sum_{n \not= 0} \hat f[n] e^{2\pi n(kv+\varphi(v))} \,\dif v \\
& = & \int_{\TT} (\varphi'(v))^2 f(kv + \varphi(v)) \,\dif v.
\end{eqnarray*}
We bound this integral by
$\|f\|_{L^\infty(\TT)} \|\varphi'\|^2_{L^2(\TT)}$.
Combining these, the expression of (\ref{c[k]-formula})
can now be bounded from above in absolute value as
\begin{equation*}
\big |c[k] \big | \leq
\frac{1}{k^2} \left(\frac{1}{\sqrt{12}} \,\|f\|_{L^2(\TT)}\|\varphi'\|_{TV}
+ \|f\|_{L^\infty(\TT)} \|\varphi'\|^2_{L^2(\TT)} \right),
\end{equation*}
concluding the proof.
\hfill $\Box$

\section*{Acknowledgements}

The authors would like to thank Ingrid Daubechies, Ron DeVore, 
\"Ozg\"ur Y\i lmaz and Yang Wang for conversations
on the topic of \sd quantization, tiling, and related issues.


\begin{thebibliography}{12}

\bibitem{Adler}
R.~L. Adler, B.~P. Kitchens, M. Martens, C.~P. Tresser, C.~W. Wu,
\newblock ``The mathematics of halftoning,''
\newblock {\em IBM J. Res. \& Dev.}, vol.~47, no.~1, Jan 2003.

\bibitem{Anastas}
D. Anastassiou,
\newblock ``Error diffusion coding for A/D conversion,''
\newblock {\em IEEE Trans. on Circuits and Systems}, vol.~36, no.~3,
pp.~1175--1186, Sept. 1989.

\bibitem{Anzai}
H. Anzai,
\newblock ``Ergodic Skew Product Transformation on the Torus,''
\newblock {\em Osaka Math. J.}, {\bf 3} (1951), pp.~83--99.

\bibitem{Bernard}
T. Bernard,
\newblock ``From \protect{$\Sigma-\Delta$} modulation to digital halftoning
of images,''
\newblock {\em Proc. IEEE Int. Conf. on Acoustics, Speech and
Signal Proc.}, May 1991, pp. 2805--2808, Toronto.

\bibitem{ingrid_democ}
R. Calderbank and I. Daubechies,
\newblock ``The Pros and Cons of Democracy,''
\newblock {\em IEEE Trans. in Inform. Theory}, vol.~48, pp.~1721--1725,
June 2002.

\bibitem{sigmadelta1}
J.~C. Candy and G.~C. Temes, Eds.,
\newblock {\em Oversampling Delta-Sigma Data Converters: Theory, Design and
  Simulation},
\newblock IEEE Press, 1992.

\bibitem{Furst}
H. Furstenberg,
\newblock ``Strict ergodicity and transformations of the torus,''
\newblock {\em Amer. J. Math.}, vol. 83, pp. 573--601, 1961.

\bibitem{GrayMulti}
W.~Chou, P.~W. Wong, and R.~M. Gray,
\newblock ``Multistage $\Sigma\Delta$ modulation,''
\newblock {\em IEEE Trans. Inform. Theory}, vol. 35, pp. 784--796, July 1989.

\bibitem{ingrid_devore}
I.~Daubechies, R.~DeVore,
\newblock  ``Reconstructing a Bandlimited Function From Very Coarsely
Quantized Data: A Family of Stable Sigma-Delta Modulators of
Arbitrary Order'', {\em Ann. of Math.}, vol. 158, no. 2, pp. 643--674,
Sept. 2003.

\bibitem{gray1}
R.~M. Gray,
\newblock ``Oversampled sigma-delta modulation,''
\newblock {\em IEEE Trans. on Comm.}, vol.~COM-35, pp.~481--489, May 1987.

\bibitem{gray2}
R.~M. Gray,
\newblock ``Spectral Analysis of Quantization Noise in a
Single-Loop Sigma-Delta Modulator with dc Input,''
\newblock {\em IEEE Trans. on Comm.}, vol.~COM-37, pp. 588--599, June 1989.

\bibitem{CSG_thesis}
C.~S. G\"unt\"urk,
\newblock {\em Harmonic Analysis of Two Problems in Signal Quantization and
 Compression},
\newblock Ph.D. thesis, Princeton University, 2000.

\bibitem{Sinan1}
C.~S. G{\" u}nt{\" u}rk
\newblock ``Approximating a Bandlimited Function Using Very Coarsely Quantized
Data: Improved Error Estimates in Sigma-Delta
Modulation'', {\em J. Amer. Math. Soc.}, 
posted on August 1, 2003, PII S 0894-0347(03)00436-3 (to appear in print).

\bibitem{SinanThao1}
C.~S. G{\" u}nt{\" u}rk and N.~T. Thao,
\newblock ``Refined Analysis of MSE in Second Order Sigma-Delta
Modulation with DC Inputs,''
submitted to IEEE Transactions on Information Theory, in revision.

\bibitem{exp_decay}
C.~S. G{\" u}nt{\" u}rk,
\newblock ``One-Bit Sigma-Delta Quantization with Exponential Accuracy,''
\newblock {\em Comm. Pure Appl. Math.},  vol. 56, pp. 1608--1630, no. 11,
2003.

\bibitem{He92}
N.~He, F.~Kuhlmann, and A.~Buzo,
\newblock ``Multi-loop $\Sigma\Delta$ quantization,''
\newblock {\em IEEE Trans. Inform. Theory}, vol. 38, pp.~1015--1028, May 1992.

\bibitem{Katok}
A. Katok, and B. Hasselblatt,
{\em Introduction to the Modern Theory of Dynamical Systems},
\newblock Cambridge University Press, 1995.

\bibitem{KatzN}
Y. Katznelson,
{\em An Introduction to Harmonic Analysis},
\newblock John Wiley \& Sons, 1968 (reprint: Dover Publs. Inc.).

\bibitem{Kite1}
T.~D. Kite, B.~L. Evans, A.~C. Bovik, and T.~L. Sculley,
\newblock ``Digital image halftoning as 2-D Delta-Sigma modulation'',
\newblock {\em Proc. IEEE Int. Conf. on
Image Proc.}, vol. I, pp. 799-802, Oct. 26-29, 1997, Santa Barbara, CA.

\bibitem{kuipers}
L.~Kuipers and H.~Niederreiter,
\newblock {\em Uniform Distribution of Sequences},
\newblock Wiley, 1974.

\bibitem{sigmadelta2}
S.~R. Norsworthy, R.~Schreier, and G.~C. Temes, Eds.,
\newblock {\em Delta-Sigma Data Converters: Theory, Design and Simulation},
\newblock IEEE Press, 1996.

\bibitem{Parry}
W.~Parry,
\newblock {\em Topics in Ergodic Theory},
\newblock Cambridge University Press, 1981.

\bibitem{Schreier}
R.~Schreier, M.~V. Goodson, and B.~Zhang,
\newblock ``An algorithm for computing convex positively invariant sets 
for delta-sigma modulators,''
\newblock {\em IEEE Trans. on Circuits and Systems, I}, vol. 44, pp.~38--44, 
Jan. 1997.

\bibitem{ThaoBreaks}
N.~T. Thao,
\newblock ``Breaking the feedback loop of a class of $\Sigma\Delta$ A/D
converters,''
\newblock {\em IEEE Trans. Signal Processing}, submitted.

\bibitem{Ulichney}
R. Ulichney,
\newblock {\em Digital Halftoning},
\newblock MIT Press, Cambridge, 1987.

\bibitem{Yilmaz}
\"O.~Y\i lmaz,
\newblock ``Stability analysis for several second-order sigma-delta methods 
of coarse quantization of bandlimited functions,''
\newblock Constr. Approx. 18 (2002), no. 4, 599--623.


\end{thebibliography}
\end{document}